\documentclass[12pt,final]{amsart}

\input epsf
\theoremstyle{plain}
\newtheorem{thm}[subsubsection]{Theorem}
\newtheorem{lem}[subsubsection]{Lemma}
\newtheorem{prop}[subsubsection]{Proposition}
\newtheorem{cor}[subsubsection]{Corollary}

\theoremstyle{definition}
\newtheorem{rem}[subsubsection]{Remark}
\newtheorem{defn}[subsubsection]{Definition}

\newcommand{\moduli}{\mathcal{M}}
\newcommand{\modulitwo}{{\mathcal{M}_2}}

\def\ui{\underline{i}}

\newcommand{\fs}{\mbox{$/\!\!/$}}

\def\L\mathcal{{L}}
\def\A\mathcal{{A}}
\def\a{\alpha}
\def\b{\beta}

\def\e{\epsilon}
\def\f{\phi}

\def\L{\Lambda}

\def\t{\tau}
\def\te{\theta}

\def\ui{\underline{i}}
\def\uj{\underline{j}}

\def\ni{\noindent}
\def\vs{\vskip}

\edef\thinlines{\the\catcode`@ }%
\catcode`@ = 11
\let\@oldatcatcode = \thinlines

\def\smash@@{\relax 
  \ifmmode\def\next{\mathpalette\mathsm@sh}\else\let\next\makesm@sh
  \fi\next}
\def\makesm@sh#1{\setbox\z@\hbox{#1}\finsm@sh}
\def\mathsm@sh#1#2{\setbox\z@\hbox{$\m@th#1{#2}$}\finsm@sh}
\def\finsm@sh{\ht\z@\z@ \dp\z@\z@ \box\z@}

\edef\@oldandcatcode{\the\catcode`& }%
\catcode`& = 11

\def\&whilenoop#1{}%
\def\&whiledim#1\do #2{\ifdim #1\relax#2\&iwhiledim{#1\relax#2}\fi}%
\def\&iwhiledim#1{\ifdim #1\let\&nextwhile=\&iwhiledim
        \else\let\&nextwhile=\&whilenoop\fi\&nextwhile{#1}}%

\newif\if&negarg
\newdimen\&wholewidth
\newdimen\&halfwidth

\font\tenln=line10

\def\thinlines{\let\&linefnt\tenln \let\&circlefnt\tencirc
  \&wholewidth\fontdimen8\tenln \&halfwidth .5\&wholewidth}%
\def\thicklines{\let\&linefnt\tenlnw \let\&circlefnt\tencircw
  \&wholewidth\fontdimen8\tenlnw \&halfwidth .5\&wholewidth}%

\def\drawline(#1,#2)#3{\&xarg #1\relax \&yarg #2\relax \&linelen=#3\relax
  \ifnum\&xarg =0 \&vline \else \ifnum\&yarg =0 \&hline \else \&sline\fi\fi}%

\def\&sline{\leavevmode
  \ifnum\&xarg< 0 \&negargtrue \&xarg -\&xarg \&yyarg -\&yarg
  \else \&negargfalse \&yyarg \&yarg \fi
  \ifnum \&yyarg >0 \&tempcnta\&yyarg \else \&tempcnta -\&yyarg \fi
  \ifnum\&tempcnta>6 \&badlinearg \&yyarg0 \fi
  \ifnum\&xarg>6 \&badlinearg \&xarg1 \fi
  \setbox\&linechar\hbox{\&linefnt\&getlinechar(\&xarg,\&yyarg)}%
  \ifnum \&yyarg >0 \let\&upordown\raise \&clnht\z@
  \else\let\&upordown\lower \&clnht \ht\&linechar\fi
  \&clnwd=\wd\&linechar
  \&whiledim \&clnwd <\&linelen \do {%
    \&upordown\&clnht\copy\&linechar
    \advance\&clnht \ht\&linechar
    \advance\&clnwd \wd\&linechar
  }%
  \advance\&clnht -\ht\&linechar
  \advance\&clnwd -\wd\&linechar
  \&tempdima\&linelen\advance\&tempdima -\&clnwd
  \&tempdimb\&tempdima\advance\&tempdimb -\wd\&linechar
  \hskip\&tempdimb \multiply\&tempdima \@m
  \&tempcnta \&tempdima \&tempdima \wd\&linechar \divide\&tempcnta \&tempdima
  \&tempdima \ht\&linechar \multiply\&tempdima \&tempcnta
  \divide\&tempdima \@m
  \advance\&clnht \&tempdima
  \ifdim \&linelen <\wd\&linechar \hskip \wd\&linechar
  \else\&upordown\&clnht\copy\&linechar\fi}%

\def\&hline{\vrule height \&halfwidth depth \&halfwidth width \&linelen}%

\def\&getlinechar(#1,#2){\&tempcnta#1\relax\multiply\&tempcnta 8
  \advance\&tempcnta -9 \ifnum #2>0 \advance\&tempcnta #2\relax\else
  \advance\&tempcnta -#2\relax\advance\&tempcnta 64 \fi
  \char\&tempcnta}%

\def\drawvector(#1,#2)#3{\&xarg #1\relax \&yarg #2\relax
  \&tempcnta \ifnum\&xarg<0 -\&xarg\else\&xarg\fi
  \ifnum\&tempcnta<5\relax \&linelen=#3\relax
    \ifnum\&xarg =0 \&vvector \else \ifnum\&yarg =0 \&hvector
    \else \&svector\fi\fi\else\&badlinearg\fi}%

\def\&hvector{\ifnum\&xarg<0 \rlap{\&linefnt\&getlarrow(1,0)}\fi \&hline
  \ifnum\&xarg>0 \llap{\&linefnt\&getrarrow(1,0)}\fi}%

\def\&vvector{\ifnum \&yarg <0 \&downvector \else \&upvector \fi}%

\def\&svector{\&sline
  \&tempcnta\&yarg \ifnum\&tempcnta <0 \&tempcnta=-\&tempcnta\fi
  \ifnum\&tempcnta <5
    \if&negarg\ifnum\&yarg>0                   
      \llap{\lower\ht\&linechar\hbox to\&linelen{\&linefnt
        \&getlarrow(\&xarg,\&yyarg)\hss}}\else 
      \llap{\hbox to\&linelen{\&linefnt\&getlarrow(\&xarg,\&yyarg)\hss}}\fi
    \else\ifnum\&yarg>0                        
      \&tempdima\&linelen \multiply\&tempdima\&yarg
      \divide\&tempdima\&xarg \advance\&tempdima-\ht\&linechar
      \raise\&tempdima\llap{\&linefnt\&getrarrow(\&xarg,\&yyarg)}\else
      \&tempdima\&linelen \multiply\&tempdima-\&yarg 
      \divide\&tempdima\&xarg
      \lower\&tempdima\llap{\&linefnt\&getrarrow(\&xarg,\&yyarg)}\fi\fi
  \else\&badlinearg\fi}%

\def\&getlarrow(#1,#2){\ifnum #2 =\z@ \&tempcnta='33\else
\&tempcnta=#1\relax\multiply\&tempcnta \sixt@@n \advance\&tempcnta
-9 \&tempcntb=#2\relax\multiply\&tempcntb \tw@ \ifnum \&tempcntb
>0 \advance\&tempcnta \&tempcntb\relax \else\advance\&tempcnta
-\&tempcntb\advance\&tempcnta 64
\fi\fi\char\&tempcnta}%

\def\&getrarrow(#1,#2){\&tempcntb=#2\relax
\ifnum\&tempcntb < 0 \&tempcntb=-\&tempcntb\relax\fi \ifcase
\&tempcntb\relax \&tempcnta='55 \or \ifnum #1<3
\&tempcnta=#1\relax\multiply\&tempcnta 24 \advance\&tempcnta -6
\else \ifnum #1=3 \&tempcnta=49 \else\&tempcnta=58 \fi\fi\or
\ifnum #1<3 \&tempcnta=#1\relax\multiply\&tempcnta 24
\advance\&tempcnta -3 \else \&tempcnta=51\fi\or
\&tempcnta=#1\relax\multiply\&tempcnta \sixt@@n \advance\&tempcnta
-\tw@ \else \&tempcnta=#1\relax\multiply\&tempcnta \sixt@@n
\advance\&tempcnta 7 \fi\ifnum #2<0 \advance\&tempcnta 64 \fi
\char\&tempcnta}%

\def\&vline{\ifnum \&yarg <0 \&downline \else \&upline\fi}%

\def\&upline{\hbox to \z@{\hskip -\&halfwidth \vrule width \&wholewidth
   height \&linelen depth \z@\hss}}%

\def\&downline{\hbox to \z@{\hskip -\&halfwidth \vrule width \&wholewidth
   height \z@ depth \&linelen \hss}}%

\def\&upvector{\&upline\setbox\&tempboxa\hbox{\&linefnt\char'66}\raise
     \&linelen \hbox to\z@{\lower \ht\&tempboxa\box\&tempboxa\hss}}%

\def\&downvector{\&downline\lower \&linelen
      \hbox to \z@{\&linefnt\char'77\hss}}%

\def\&badlinearg{\errmessage{Bad \string\arrow\space argument.}}%

\thinlines

\countdef\&xarg     0 \countdef\&yarg     2 \countdef\&yyarg    4
\countdef\&tempcnta 6 \countdef\&tempcntb 8

\dimendef\&linelen  0 \dimendef\&clnwd    2 \dimendef\&clnht    4
\dimendef\&tempdima 6 \dimendef\&tempdimb 8

\chardef\@arrbox    0 \chardef\&linechar  2
\chardef\&tempboxa  2           

\let\lft^%
\let\rt_

\newif\if@pslope 
\def\@findslope(#1,#2){\ifnum#1>0
  \ifnum#2>0 \@pslopetrue \else\@pslopefalse\fi \else
  \ifnum#2>0 \@pslopefalse \else\@pslopetrue\fi\fi}%

\def\generalsmap(#1,#2){\getm@rphposn(#1,#2)\plnmorph\futurelet\next\addm@rph}%

\def\sline(#1,#2){\setbox\@arrbox=\hbox{\drawline(#1,#2){\sarrowlength}}%
  \@findslope(#1,#2)\d@@blearrfalse\generalsmap(#1,#2)}%
\def\arrow(#1,#2){\setbox\@arrbox=\hbox{\drawvector(#1,#2){\sarrowlength}}%
  \@findslope(#1,#2)\d@@blearrfalse\generalsmap(#1,#2)}%

\newif\ifd@@blearr

\def\bisline(#1,#2){\@findslope(#1,#2)%
  \if@pslope \let\@upordown\raise \else \let\@upordown\lower\fi
  \getch@nnel(#1,#2)\setbox\@arrbox=\hbox{\@upordown\@vchannel
    \rlap{\drawline(#1,#2){\sarrowlength}}%
      \hskip\@hchannel\hbox{\drawline(#1,#2){\sarrowlength}}}%
  \d@@blearrtrue\generalsmap(#1,#2)}%
\def\biarrow(#1,#2){\@findslope(#1,#2)%
  \if@pslope \let\@upordown\raise \else \let\@upordown\lower\fi
  \getch@nnel(#1,#2)\setbox\@arrbox=\hbox{\@upordown\@vchannel
    \rlap{\drawvector(#1,#2){\sarrowlength}}%
      \hskip\@hchannel\hbox{\drawvector(#1,#2){\sarrowlength}}}%
  \d@@blearrtrue\generalsmap(#1,#2)}%
\def\adjarrow(#1,#2){\@findslope(#1,#2)%
  \if@pslope \let\@upordown\raise \else \let\@upordown\lower\fi
  \getch@nnel(#1,#2)\setbox\@arrbox=\hbox{\@upordown\@vchannel
    \rlap{\drawvector(#1,#2){\sarrowlength}}%
      \hskip\@hchannel\hbox{\drawvector(-#1,-#2){\sarrowlength}}}%
  \d@@blearrtrue\generalsmap(#1,#2)}%

\newif\ifrtm@rph
\def\@shiftmorph#1{\hbox{\setbox0=\hbox{$\scriptstyle#1$}%
  \setbox1=\hbox{\hskip\@hm@rphshift\raise\@vm@rphshift\copy0}%
  \wd1=\wd0 \ht1=\ht0 \dp1=\dp0 \box1}}%
\def\@hm@rphshift{\ifrtm@rph
  \ifdim\hmorphposnrt=\z@\hmorphposn\else\hmorphposnrt\fi \else
  \ifdim\hmorphposnlft=\z@\hmorphposn\else\hmorphposnlft\fi \fi}%
\def\@vm@rphshift{\ifrtm@rph
  \ifdim\vmorphposnrt=\z@\vmorphposn\else\vmorphposnrt\fi \else
  \ifdim\vmorphposnlft=\z@\vmorphposn\else\vmorphposnlft\fi \fi}%

\def\addm@rph{\ifx\next\lft\let\temp=\lftmorph\else
  \ifx\next\rt\let\temp=\rtmorph\else\let\temp\relax\fi\fi \temp}%

\def\plnmorph{\dimen1\wd\@arrbox \ifdim\dimen1<\z@ \dimen1-\dimen1\fi
  \vcenter{\box\@arrbox}}%
\def\lftmorph\lft#1{\rtm@rphfalse \setbox0=\@shiftmorph{#1}%
  \if@pslope \let\@upordown\raise \else \let\@upordown\lower\fi
  \llap{\@upordown\@vmorphdflt\hbox to\dimen1{\hss 
    \llap{\box0}\hss}\hskip\@hmorphdflt}\futurelet\next\addm@rph}%
\def\rtmorph\rt#1{\rtm@rphtrue \setbox0=\@shiftmorph{#1}%
  \if@pslope \let\@upordown\lower \else \let\@upordown\raise\fi
  \llap{\@upordown\@vmorphdflt\hbox to\dimen1{\hss
    \rlap{\box0}\hss}\hskip-\@hmorphdflt}\futurelet\next\addm@rph}%

\def\getm@rphposn(#1,#2){\ifd@@blearr \dimen@\morphdist \advance\dimen@ by
  .5\channelwidth \@getshift(#1,#2){\@hmorphdflt}{\@vmorphdflt}{\dimen@}\else
  \@getshift(#1,#2){\@hmorphdflt}{\@vmorphdflt}{\morphdist}\fi}%

\def\getch@nnel(#1,#2){\ifdim\hchannel=\z@ \ifdim\vchannel=\z@
    \@getshift(#1,#2){\@hchannel}{\@vchannel}{\channelwidth}%
    \else \@hchannel\hchannel \@vchannel\vchannel \fi
  \else \@hchannel\hchannel \@vchannel\vchannel \fi}%

\def\@getshift(#1,#2)#3#4#5{\dimen@ #5\relax
  \&xarg #1\relax \&yarg #2\relax
  \ifnum\&xarg<0 \&xarg -\&xarg \fi
  \ifnum\&yarg<0 \&yarg -\&yarg \fi
  \ifnum\&xarg<\&yarg \&negargtrue \&yyarg\&xarg \&xarg\&yarg \&yarg\&yyarg\fi
  \ifcase\&xarg \or  
    \ifcase\&yarg    
      \dimen@i \z@ \dimen@ii \dimen@ \or 
      \dimen@i .7071\dimen@ \dimen@ii .7071\dimen@ \fi \or
    \ifcase\&yarg    
      \or 
      \dimen@i .4472\dimen@ \dimen@ii .8944\dimen@ \fi \or
    \ifcase\&yarg    
      \or 
      \dimen@i .3162\dimen@ \dimen@ii .9486\dimen@ \or
      \dimen@i .5547\dimen@ \dimen@ii .8321\dimen@ \fi \or
    \ifcase\&yarg    
      \or 
      \dimen@i .2425\dimen@ \dimen@ii .9701\dimen@ \or\or
      \dimen@i .6\dimen@ \dimen@ii .8\dimen@ \fi \or
    \ifcase\&yarg    
      \or 
      \dimen@i .1961\dimen@ \dimen@ii .9801\dimen@ \or
      \dimen@i .3714\dimen@ \dimen@ii .9284\dimen@ \or
      \dimen@i .5144\dimen@ \dimen@ii .8575\dimen@ \or
      \dimen@i .6247\dimen@ \dimen@ii .7801\dimen@ \fi \or
    \ifcase\&yarg    
      \or 
      \dimen@i .1645\dimen@ \dimen@ii .9864\dimen@ \or\or\or\or
      \dimen@i .6402\dimen@ \dimen@ii .7682\dimen@ \fi \fi
  \if&negarg \&tempdima\dimen@i \dimen@i\dimen@ii \dimen@ii\&tempdima\fi
  #3\dimen@i\relax #4\dimen@ii\relax }%

\catcode`\&=4  

\def\generalhmap{\futurelet\next\@generalhmap}%
\def\@generalhmap{\ifx\next^ \let\temp\generalhm@rph\else
  \ifx\next_ \let\temp\generalhm@rph\else \let\temp\m@kehmap\fi\fi \temp}%
\def\generalhm@rph#1#2{\ifx#1^
    \toks@=\expandafter{\the\toks@#1{\rtm@rphtrue\@shiftmorph{#2}}}\else
    \toks@=\expandafter{\the\toks@#1{\rtm@rphfalse\@shiftmorph{#2}}}\fi
  \generalhmap}%
\def\m@kehmap{\mathrel{\smash@@{\the\toks@}}}%

\def\mapright{\toks@={\mathop{\vcenter{\smash@@{\drawrightarrow}}}\limits}%
  \generalhmap}%
\def\mapleft{\toks@={\mathop{\vcenter{\smash@@{\drawleftarrow}}}\limits}%
  \generalhmap}%
\def\bimapright{\toks@={\mathop{\vcenter{\smash@@{\drawbirightarrow}}}\limits}%
  \generalhmap}%
\def\bimapleft{\toks@={\mathop{\vcenter{\smash@@{\drawbileftarrow}}}\limits}%
  \generalhmap}%
\def\adjmapright{\toks@={\mathop{\vcenter{\smash@@{\drawadjrightarrow}}}\limits}%
  \generalhmap}%
\def\adjmapleft{\toks@={\mathop{\vcenter{\smash@@{\drawadjleftarrow}}}\limits}%
  \generalhmap}%
\def\hline{\toks@={\mathop{\vcenter{\smash@@{\drawhline}}}\limits}%
  \generalhmap}%
\def\bihline{\toks@={\mathop{\vcenter{\smash@@{\drawbihline}}}\limits}%
  \generalhmap}%

\def\drawrightarrow{\hbox{\drawvector(1,0){\harrowlength}}}%
\def\drawleftarrow{\hbox{\drawvector(-1,0){\harrowlength}}}%
\def\drawbirightarrow{\hbox{\raise.5\channelwidth
  \hbox{\drawvector(1,0){\harrowlength}}\lower.5\channelwidth
  \llap{\drawvector(1,0){\harrowlength}}}}%
\def\drawbileftarrow{\hbox{\raise.5\channelwidth
  \hbox{\drawvector(-1,0){\harrowlength}}\lower.5\channelwidth
  \llap{\drawvector(-1,0){\harrowlength}}}}%
\def\drawadjrightarrow{\hbox{\raise.5\channelwidth
  \hbox{\drawvector(-1,0){\harrowlength}}\lower.5\channelwidth
  \llap{\drawvector(1,0){\harrowlength}}}}%
\def\drawadjleftarrow{\hbox{\raise.5\channelwidth
  \hbox{\drawvector(1,0){\harrowlength}}\lower.5\channelwidth
  \llap{\drawvector(-1,0){\harrowlength}}}}%
\def\drawhline{\hbox{\drawline(1,0){\harrowlength}}}%
\def\drawbihline{\hbox{\raise.5\channelwidth
  \hbox{\drawline(1,0){\harrowlength}}\lower.5\channelwidth
  \llap{\drawline(1,0){\harrowlength}}}}%

\def\generalvmap{\futurelet\next\@generalvmap}%
\def\@generalvmap{\ifx\next\lft \let\temp\generalvm@rph\else
  \ifx\next\rt \let\temp\generalvm@rph\else \let\temp\m@kevmap\fi\fi \temp}%
\toksdef\toks@@=1
\def\generalvm@rph#1#2{\ifx#1\rt 
    \toks@=\expandafter{\the\toks@
      \rlap{$\vcenter{\rtm@rphtrue\@shiftmorph{#2}}$}}\else 
    \toks@@={\llap{$\vcenter{\rtm@rphfalse\@shiftmorph{#2}}$}}%
    \toks@=\expandafter\expandafter\expandafter{\expandafter\the\expandafter
      \toks@@ \the\toks@}\fi \generalvmap}%
\def\m@kevmap{\the\toks@}%

\def\mapdown{\toks@={\vcenter{\drawdownarrow}}\generalvmap}%
\def\mapup{\toks@={\vcenter{\drawuparrow}}\generalvmap}%
\def\bimapdown{\toks@={\vcenter{\drawbidownarrow}}\generalvmap}%
\def\bimapup{\toks@={\vcenter{\drawbiuparrow}}\generalvmap}%
\def\adjmapdown{\toks@={\vcenter{\drawadjdownarrow}}\generalvmap}%
\def\adjmapup{\toks@={\vcenter{\drawadjuparrow}}\generalvmap}%
\def\vline{\toks@={\vcenter{\drawvline}}\generalvmap}%
\def\bivline{\toks@={\vcenter{\drawbivline}}\generalvmap}%

\def\drawdownarrow{\hbox to5pt{\hss\drawvector(0,-1){\varrowlength}\hss}}%
\def\drawuparrow{\hbox to5pt{\hss\drawvector(0,1){\varrowlength}\hss}}%
\def\drawbidownarrow{\hbox to5pt{\hss\hbox{\drawvector(0,-1){\varrowlength}}%
  \hskip\channelwidth\hbox{\drawvector(0,-1){\varrowlength}}\hss}}%
\def\drawbiuparrow{\hbox to5pt{\hss\hbox{\drawvector(0,1){\varrowlength}}%
  \hskip\channelwidth\hbox{\drawvector(0,1){\varrowlength}}\hss}}%
\def\drawadjdownarrow{\hbox to5pt{\hss\hbox{\drawvector(0,-1){\varrowlength}}%
  \hskip\channelwidth\lower\varrowlength
  \hbox{\drawvector(0,1){\varrowlength}}\hss}}%
\def\drawadjuparrow{\hbox to5pt{\hss\hbox{\drawvector(0,1){\varrowlength}}%
  \hskip\channelwidth\raise\varrowlength
  \hbox{\drawvector(0,-1){\varrowlength}}\hss}}%
\def\drawvline{\hbox to5pt{\hss\drawline(0,1){\varrowlength}\hss}}%
\def\drawbivline{\hbox to5pt{\hss\hbox{\drawline(0,1){\varrowlength}}%
  \hskip\channelwidth\hbox{\drawline(0,1){\varrowlength}}\hss}}%

\def\commdiag#1{\null\,
  \vcenter{\commdiagbaselines
  \m@th\ialign{\hfil$##$\hfil&&\hfil$\mkern4mu ##$\hfil\crcr
      \mathstrut\crcr\noalign{\kern-\baselineskip}
      #1\crcr\mathstrut\crcr\noalign{\kern-\baselineskip}}}\,}%

\def\commdiagbaselines{\baselineskip15pt \lineskip3pt \lineskiplimit3pt }%
\def\gridcommdiag#1{\null\,
  \vcenter{\offinterlineskip
  \m@th\ialign{&\vbox to\vgrid{\vss
    \hbox to\hgrid{\hss\smash@@{$##$}\hss}}\crcr
      \mathstrut\crcr\noalign{\kern-\vgrid}
      #1\crcr\mathstrut\crcr\noalign{\kern-.5\vgrid}}}\,}%

\newdimen\harrowlength \harrowlength=60pt
\newdimen\varrowlength \varrowlength=.618\harrowlength
\newdimen\sarrowlength \sarrowlength=\harrowlength

\newdimen\hmorphposn \hmorphposn=\z@
\newdimen\vmorphposn \vmorphposn=\z@
\newdimen\morphdist  \morphdist=4pt

\dimendef\@hmorphdflt 0       
\dimendef\@vmorphdflt 2       

\newdimen\hmorphposnrt  \hmorphposnrt=\z@
\newdimen\hmorphposnlft \hmorphposnlft=\z@
\newdimen\vmorphposnrt  \vmorphposnrt=\z@
\newdimen\vmorphposnlft \vmorphposnlft=\z@

\newdimen\hgrid \hgrid=15pt
\newdimen\vgrid \vgrid=15pt

\newdimen\hchannel  \hchannel=0pt
\newdimen\vchannel  \vchannel=0pt
\newdimen\channelwidth \channelwidth=3pt

\dimendef\@hchannel 0         
\dimendef\@vchannel 2         

\catcode`& = \@oldandcatcode \catcode`@ = \@oldatcatcode

\input epsf
\begin{document}

\title[Standard monomial bases]{Standard monomial bases, Moduli spaces of vector bundles
 \& Invariant theory}
\author[V. Lakshmibai]{V. Lakshmibai${}^{\dag}$}
\address{Department of Mathematics\\ Northeastern University\\
Boston, MA 02115} \email{lakshmibai@neu.edu}
\thanks{${}^{\dag}$ Partially suported
by NSF grant DMS-0400679 and NSA-MDA904-03-1-0034.}
\author{K.N. Raghavan}
\address{The Institute of Mathematical Sciences, Chennai 600\,113, INDIA} \email{knr@imsc.res.in}
\author{P. Sankaran}
\address{The Institute of Mathematical Sciences, Chennai 600\,113, INDIA} \email{sankaran@imsc.res.in}
\author{P. Shukla}
\address{Department of Mathematics\\
Suffolk University\\ Boston, MA 02114}
\email{shukla@mcs.suffolk.edu}

\begin{abstract}
 Consider the diagonal action of $SO_n(K)$ on the affine space $X=V^{\oplus m}$ where $V=K^n,\,K$
 an algebraically closed field of characteristic $\not= 2$. We construct a
 ``standard monomial" basis for the ring of invariants
 $K[X]^{SO_n(K)}$. As a consequence, we deduce that
 $K[X]^{SO_n(K)}$ is Cohen-Macaulay.
 As the first application, we present the first and second fundamental theorems for
$SO_n(K)$-actions. As the second application, assuming that the
characteristic of $K$ is $\neq 2,3$, we give a
characteristic-free proof of the Cohen-Macaulayness of the moduli
space $\mathcal{M}_2$ of equivalence classes of semi-stable, rank
$2$, degree $0$ vector bundles on a smooth projective curve of
genus $> 2$. As the third application, we describe a $K$-basis for
the ring of invariants for the adjoint action of $SL_2(K)$ on $m$
copies of $sl_2(K)$ in terms of traces.
\end{abstract}

\maketitle


\section*{Introduction} This paper is a sequel to \cite{l-s}. Let $V=K^n$, together with a symmetric bilinear form
$\langle,\rangle$, $K$ being an algebraically closed field of
characteristic $\not= 2$. Let  $X=\underset{m\text{
copies}}{\underbrace{V\oplus \cdots \oplus V}}$. In \cite{d-p},
a characteristic-free basis is described for
$K[X]^{O_n(K)}$ (for the diagonal action of $O_n(K)$ on $X$). 
The diagonal action of $SO_n(K)$ on $X$ is also considered, and
a set of algebra generators is described for $K[X]^{SO_n(K)}$ (cf.
\cite{d-p}, Theorem 5.6,(2)).

The main goal of this paper is to prove the Cohen-Macaulayness of
$K[X]^{SO_n(K)}$ (note that the Cohen-Macaulayness of
$K[X]^{O_n(K)}$ follows from the fact that

\ni $Spec\,(K[X]^{O_n(K)})$ is a certain determinantal variety
inside Sym$\,M_{m}$, the space of symmetric $m\times m$ matrices;
note also that in characteristic $0$, the
 Cohen-Macaulayness of $K[X]^{SO_n(K)}$ follows from
\cite{h-r,bo}). We adopt the ``deformation technique" for this
purpose. As mentioned in the introduction of \cite{l-s}, the
``deformation technique" has proven to be quite effective for
proving the Cohen-Macaulayness of algebraic varieties. This
technique consists in constructing a flat family over
$\mathbb{A}^1$, with the given variety as the generic fiber
(corresponding to $t\in K$ invertible). If the special fiber
(corresponding to $t=0$) is Cohen-Macaulay, then one may conclude
the Cohen-Macaulayness of the given variety. Hodge algebras (cf.
\cite{d-e-p}) are typical examples where the deformation technique
affords itself very well. Deformation technique is also used in
\cite{d-l,h-l,g-l,chi,cald,l-s}.
 The philosophy behind these works is that if there is a
``standard monomial basis" for the co-ordinate ring of the given
variety, then the deformation technique will work well in general
(using the ``straightening relations").
It is this philosophy that we adopt in this paper in proving the
Cohen-Macaulayness of $K[X]^{SO_n(K)}$. To be more precise, the
proof of the Cohen-Macaulayness of $K[X]^{SO_n(K)}$ is
accomplished in the following steps:

\vs.2cm $\bullet$\hskip1cm We first construct a $K$-subalgebra $S$
of $K[X]^{SO_n(K)}$ by prescribing a set of algebra generators
$\{f_\alpha,\alpha\in D\}, D$ being a doset associated to a finite
partially ordered set $P$ (here, ``doset" is as defined in
\cite{d-l}; see also Definition \ref{dos}).

$\bullet$\hskip1cm We construct a ``standard monomial" basis for
$S$ by

(i) defining ``standard monomials" in the $f_\alpha$'s (cf.
Definition \ref{std})

(ii) writing down the straightening relation for a non-standard
(degree $2$) monomial $f_\alpha f_\beta $ (cf. Proposition
\ref{relations})

(iii) proving linear independence of standard monomials (by
relating the generators of $S$ to certain determinantal varieties
inside the space of symmetric matrices) (cf. Proposition
\ref{indpt})

(iv) proving the generation of $S$ (as a vector space) by standard
monomials (using (ii)). In fact, to prove the generation for $S$,
we first prove generation for a ``graded version" $R(D)$ of $S$,
where $D$ as above is a doset associated to a  finite partially
ordered set $P$. We then deduce the generation for $S$. In fact,
we construct a ``standard monomial" basis for $R(D)$. While the
generation by standard monomials for $S$ is deduced from the
generation by standard monomials for $R(D)$, the linear
independence of standard monomials in $R(D)$ is deduced from the
linear independence of standard monomials in $S$ (cf. (iii)
above).

$\bullet$\hskip1cm We give a presentation for $S$ as a $K$-algebra
(cf. Theorem \ref{present'}).

$\bullet$\hskip1cm We prove the Cohen-Macaulayness (cf.
Proposition \ref{dose}, Corollary \ref{cohen'}) of $R(D)$ by
realizing it as a ``doset algebra with straightening law" (cf.
\cite{d-l}), and using the results of \cite{d-l}.

$\bullet$\hskip1cm We deduce the Cohen-Macaulayness of $S$ from
that of $R(D)$ (cf. Theorem \ref{main}).

$\bullet$\hskip1cm We prove (cf. Proposition \ref{regular}) that
Spec$\,S$ is regular in codimension $1$ by using the fact that
$Spec\,K[X]^{SO_n(K)}\rightarrow Spec\,K[X]^{O_n(K)}$ is a
$2$-sheeted cover.

$\bullet$\hskip1cm We deduce (cf. Proposition \ref{normal'}) the
normality of $Spec\,S$ (using Serre's criterion for normality -
$Spec\,A$ is normal if and only if $A$ has $S_2$ and $R_1$).

$\bullet$\hskip1cm We then show that the inclusion $S\subseteq
K[X]^{SO_n(K)}$ is in fact an equality by showing that the
morphism $Spec\,K[X]^{SO_n(K)}\rightarrow Spec\,S$ (induced by the
inclusion $S\subseteq K[X]^{SO_n(K)}$) satisfies the hypotheses in
Zariski Main Theorem (and hence is an isomorphism (cf. Theorem
\ref{equal})). We also deduce the the Cohen-Macaulayness of
$K[X]^{SO_n(K)}$ (cf. Theorem \ref{main'})

As the first set of main consequences, we present

$\bullet$\hskip1cm \textbf{First fundamental Theorem for
$SO_n(K)$-invariants}, i.e., describing algebra generators for
$K[X]^{SO_n(K)}$.

$\bullet$\hskip1cm \textbf{Second fundamental Theorem for
$SO_n(K)$-invariants}, i.e., describing  generators for the ideal
of relations among these algebra generators for  $K[X]^{SO_n(K)}$.

$\bullet$\hskip1cm \textbf{A standard monomial basis for
$K[X]^{SO_n(K)}$}.

As the second main consequence,
assuming that the characteristic of the
base field is $\neq 2,3$, 
we give (cf. \S \ref{moduli}) a
characteristic-free proof,
of the Cohen-Macaulayness of the moduli
space $\mathcal{M}_2$ of equivalence classes of semi-stable rank
$2$, degree $0$ vector bundles on a smooth projective curve of
genus $> 2$ by relating it to $K[X]^{SO_3(K)}$. In \cite{m-r1},
the Cohen-Macaulayness for $\mathcal{M}_2$ is deduced by proving the
Frobenius-split properties for $\mathcal{M}_2$.

As the third main consequence, we describe (cf. Theorem
\ref{standard}) a characteristic-free basis for the ring of
invariants for the (diagonal) adjoint action of $SL_2(K)$ on
$\underset{m\text{ copies}}{\underbrace{sl_2(K)\oplus \cdots
\oplus sl_2(K)}}$.

Our main goal in this paper is to prove the Cohen-Macaulayness of
$K[X]^{SO_n(K)}$; as mentioned above, this is accomplished by
first constructing a ``standard monomial" basis for the subalgebra
$S$ of $K[X]^{SO_n(K)}$, deducing Cohen-Macaulayness of $S$, and
then proving that $S$ in fact equals $K[X]^{SO_n(K)}$. Thus we
${\underline{\mathrm{do\ not}}}$ use the results of \cite{d-p}
(especially, Theorem 5.6 of \cite{d-p}), we rather give a
different proof of Theorem 5.6 of \cite{d-p}. Further, using Lemma
\ref{normality}, we get a GIT-theoretic proof of the first and
second fundamental theorems for the $O_n(K)$-action in arbitrary
characteristics which we have included in \S \ref{son}. (For the
discussions in \S \ref{alg} we need the results on the ring of
invariants for the $O_n(K)$-action - specifically, first and
second fundamental theorems for the $O_n(K)$-action.)

The sections are organized as follows. In \S \ref{prilm}, after
recalling some results (pertaining to standard monomial basis) for
Schubert varieties in the Lagrangian Grassmannian and symmetric
determinantal varieties (i.e, determinantal varieties inside the
space of symmetric matrices), we derive the straightening
relations for certain degree $2$ non-standard monomials. In \S
\ref{son}, we present a GIT-theoretic proof of of the first and
second fundamental theorems for the $O_n(K)$-action in arbitrary
characteristics. In \S \ref{alg}, we introduce the algebra $S$
($\subseteq R^{SO_n(K)}$) by describing the algebra generators,
define standard monomials in the algebra generators, and prove the
linear independence of standard monomials.  In \S \ref{alg'}, we
 introduce the algebra $R(D)$, construct a standard
monomial basis for $R(D)$, and deduce that standard monomials in
$S$ give a vector-space basis for $S$. In \S \ref{cohen}, we first
prove the Cohen-Macaulayness of $R(D)$, and then deduce the
Cohen-Macaulayness of $S$. In \S \ref{normal}, we first prove that
$S$ is normal, then show that the inclusion-induced morphism
$Spec\,R^{SO_n(K)}\longrightarrow Spec\,S$ satisfies the hypotheses
in Zariski Main Theorem, and then deduce that the inclusion
$S\subseteq R^{SO_n(K)}$ is an equality, i.e., $R^{SO_n(K)}=S$. In
\S \ref{moduli}, we present the results for the moduli space of
rank $2$ vector bundles on a smooth projective curve. In \S
\ref{sl2}, we give a characteristic-free basis for the ring of
invariants for the (diagonal) adjoint action of $SL_2(K)$ on
$\underset{m\text{ copies}}{\underbrace{sl_2(K)\oplus \cdots
\oplus sl_2(K)}}$ in terms of monomials in the traces.

\vs.2cm We thank C.S. Seshadri for many useful discussions,
especially for the discussion in~\S\ref{moduli}.    We also thank
the referees for some useful comments.   Part of this work was carried
out while the first author was visiting Chennai Mathematical Institute.
The first author wishes to express her thanks to Chennai Mathematical
Institute for the hospitality shown to her during her visit.

\section{Preliminaries}\label{prilm} In this section, we recollect some basic
results on symmetric determinantal varieties (i.e., determinantal
varieties inside Sym$\,M_m(K)$, the space of symmetric $m\times m$
matrices); specifically the standard monomial basis for the
co-ordinate rings of  symmetric determinantal varieties. Since the
results of \S \ref{stan} rely on an explicit description of the
straightening relations (of a degree $2$ non-standard monomial),
in this section we derive such straightening relations (cf.
Proposition \ref{qualitative2}) by relating symmetric
determinantal varieties to Schubert varieties in the Lagrangian
Grassmannian. We first recall some results on Schubert varieties
in the Lagrangian Grassmannian, mainly the standard monomial basis
for the homogeneous co-ordinate rings of Schubert varieties in the
Lagrangian Grassmannian (cf. \cite{g/p-2}). We then recall results
for symmetric determinantal varieties (by identifying them as open
subsets of suitable Schubert varieties in suitable Lagrangian
Grassmannians). We then derive the desired straightening
relations.

\subsection{The Lagrangian Grassmannian Variety:}\label{lag}
Let $V=K^{2m}$, ($K$ being
 the base field which we suppose to be algebraically closed of
 characteristic $\not= 2$) together with a non-degenerate,
 skew-symmetric bilinear form $\langle,\rangle$. Let $G=Sp(V)$ (the group of linear
 automorphisms of $V$ preserving $\langle,\rangle$). If $J$ is the matrix of the form,
 then $G$ may be identified with the fixed point set of the involution
 $\sigma:SL(V)\rightarrow SL(V), \sigma(A)=J^{-1}(^{t}A)^{-1}J$. Taking $J$ to be
 $$J=\begin{pmatrix}
0&J_m\\ -J_m&0
 \end{pmatrix}$$ where $J_m$ is the $m\times m$ matrix with $1$'s
 along the anti-diagonal, and $0$'s off the anti-diagonal, $B$ (resp. $T$),
 the set of upper triangular (resp. diagonal) matrices in $G$ is a
 Borel sub group (resp. a maximal torus) in $G$ (cf. \cite{stein}). Let $L_m =
 \{\mathrm{maximal\ totally\ isotropic\ sub\ spaces\ of\ } V\}$,
 \emph{the Lagrangian Grassmannian Variety}. We have a canonical
 inclusion $L_m\hookrightarrow G_{m,2m}$, where $G_{m,2m}$ is the
 Grassmannian variety of $m$-dimensional subspaces of $K^{2m}$.

 \subsection{Schubert varieties in $L_m$:}\label{jm}
 Let $I(m,2m)=\{\ui=(i_1,\dots,i_m)|1\le i_1<\dots<i_m\le 2m\}$. For $j,1\le j\le 2m$,
 let $j'=2m+1-j$. Let $$I_G(m,2m)=
 \{\ui\in I(m,2m)|\mathrm{\, for\, every\,  }j,1\le j\le 2m,
 \mathrm{\, precisely\, one\, of\, }\{j,j'\}\mathrm{\,occurs \, in\,}\ui\}$$
 Recall (cf.
 \cite{g/p-2,g/p-7}) that the Schubert varieties in $L_m$ are indexed by
 $I_G(m,2m)$; further, the partial order on the set of Schubert varieties in
 $L_m$ (given by inclusion) induces a partial order on $I_G(m,2m)$:
 $\ui\ge \uj \Leftrightarrow i_t\ge j_t , \forall t$.
  Let $w\in I_G(m,2m)$, and let $X(w)$ be the associated Schubert variety. For the
 projective embedding
 $f_m:L_m\hookrightarrow \mathbb{P}(\wedge^mV) $ (induced by the Pl\" ucker
 embedding $G_{m,2m}\hookrightarrow \mathbb{P}(\wedge^mV) $), let $R(w)$ denote
 the homogeneous co-ordinate ring of $X(w)$.

 \vs.1cm\ni\textbf{A standard monomial basis for $R(w)$:} We have
 (cf. \cite{g/p-4,g/p-5,g/p-7}; see also \cite{de}) a basis $\{p_{\tau,\varphi}\}$ for
 $R(w)_1$ indexed by admissible pairs - certain pairs $(\tau,\varphi), \tau\ge\varphi$
 of elements of $I_G(m,2m)$ (see \cite{g/p-4,g/p-5,g/p-7}
 for the definition of admissible pairs, and the description of
 $\{p_{\tau,\varphi}\}$). This basis includes the extremal weight vectors $p_\t,\tau\in I_G(m,2m)$
 corresponding to the admissible pair $(\tau,\tau)$. Thus denoting by ${\mathcal{A}}$ the set of all
admissible pairs, we have that ${\mathcal{A}}$ includes the
diagonal of $I_G(m,2m)\times I_G(m,2m)$. An admissible pair
 $(\tau,\varphi)$ such that $w\ge \tau$ is called
 \emph{an admissible pair on }$X(w)$.

 \begin{defn} A monomial of the form
$$p_{\tau_{1},\varphi_{1}}p_{\tau_{2},\varphi_{2}}\cdots
 p_{\tau_{r},\varphi_{r}},\tau_{1}\ge\varphi_{1}\ge \tau_{2}\ge \cdots \ge
 \varphi_{r}$$ is
called a
 \emph{standard monomial}.
Such a monomial is said to be \emph{standard  on } $X(w)$, if in
addition $w\ge \tau_{1}$.
\end{defn}

 Recall (cf. \cite{g/p-4,g/p-5,g/p-7})
 \begin{thm}\label{basis}
Standard monomials on $X(w)$ of degree $r$ form a basis of
$R(w)_r$.
 \end{thm}

 As a consequence, we have a qualitative description of a typical
quadratic relation on a Schubert variety $X(w)$ as given by the
following Proposition.

 \begin{prop}\label{qualitative}(\cite{d-l,g/p-4,g/p-5})
Let $(\t_1,\f_1), (\t_2,\f_2) $ be two admissible pairs on $X(w)$
such that $p_{\tau_{1},\varphi_{1}}p_{\tau_{2},\varphi_{2}}$ is
non-standard so that $\f_1\not\ge \t_2, \f_2\not\ge \t_1$. Let
\begin{equation*}
p_{\tau_{1},\varphi_{1}}p_{\tau_{2},\varphi_{2}}=\sum_{i}\ c_{i}
p_{\a_{i},\b_{i}} p_{\gamma_{i},\delta_{i}},\ c_{i}\in K^*\,
\tag{*}
\end{equation*}
be the expression for
$p_{\tau_{1},\varphi_{1}}p_{\tau_{2},\varphi_{2}}$ as a sum of
standard monomials on $X(w)$.
\begin{enumerate}
\item For every $i$, we have, $\a_i\ge$ both $\tau_{1}$ and
$\tau_{2}$. Further, for some $i$,
 if $\a_i=\tau_{1}$ (resp. $\tau_{2}$), then $\b_i> \varphi_{1}$
 (resp. $\varphi_{2}$).
\item  For every $i$, we have, $\delta_i\le$ both $\varphi_{1}$
and $\varphi_{2}$. Further, for some $i$,
 if $\delta_i=\varphi_{1}$ (resp. $\varphi_{2}$), then $\gamma_i<\tau_{1}$
 (resp. $\tau_{2}$).
 \item Suppose there exists a permutation $\sigma$ of
 $\{\tau_1,\varphi_1,\tau_2,\varphi_2\}$ such that
  $\sigma(\tau_1)\ge\sigma(\varphi_1)\ge\sigma(\tau_2)\ge\sigma(\varphi_2)
  $, then
  $(\sigma(\tau_1),\sigma(\varphi_1)),(\sigma(\tau_2),\sigma(\varphi_2))$ are
  both admissible pairs, and
  $p_{\sigma(\tau_1),\sigma(\varphi_1)}p_{\sigma(\tau_2),\sigma(\varphi_2)}$
  occurs with coefficient $\pm 1$ in (*).
 \end{enumerate}

\end{prop}

For a proof of (1), (2), we refer the reader to \cite{g/p-4},
Lemma 7.1, and of (3) to \cite{d-l}, Theorem 4.1.

\vs.2cm\ni We shall refer to a relation as in (*) as {\em a
straightening relation}.

 \vs.2cm\ni{\bf{A
presentation for $R(w)$:}} For $w\in I_G(m,2m)$, let

\ni $\mathcal{A}_w=\{(\tau,\varphi)\in \mathcal{A}\,|\,w\ge
\tau\}$. Consider the polynomial algebra
$K[x_{\tau,\varphi},\,(\tau,\varphi)\in \mathcal{A}_w]$. For two
admissible pairs $(\t_1,\f_1), (\t_2,\f_2) $ in  $\mathcal{A}_w$
such that $p_{\tau_{1},\varphi_{1}}p_{\tau_{2},\varphi_{2}}$ is
non-standard, denote $F_{(\t_1,\f_1), (\t_2,\f_2)}=x_{\t_1,\f_1},
x_{\t_2,\f_2}- \sum_{i}\ c_{i}x_{\a_{i},\b_{i}}
x_{\gamma_{i},\delta_{i}}$,

\ni $\a_i,\b_i,\gamma_{i},\delta_{i},c_i$ being as in Proposition
\ref{qualitative}. Let $I_w$ be the ideal in
$K[x_{\tau,\varphi},\,(\tau,\varphi)\in \mathcal{A}_w]$ generated
by such $F_{(\t_1,\f_1), (\t_2,\f_2)}$'s. Consider the surjective
map $f_w:K[x_{\tau,\varphi},\,(\tau,\varphi)\in
\mathcal{A}_w]\rightarrow R(w), x_{\tau,\varphi}\mapsto
p_{\tau,\varphi}$. We have

\begin{prop}\label{present}(\cite{g/p-4,g/p-5})
$ f_w$ induces an isomorphism
$K[x_{\tau,\varphi},(\tau,\varphi)\in \mathcal{A}_w]/I_w\cong
R(w)$.
\end{prop}

\subsection{The opposite big cell in
$L_m$}\label{12.1.6} We first recall the following (cf.
\cite{g/p-2}):

 \vs.2cm\ni\textbf{Fact 1.}  We have a natural embedding
 $$L_m\hookrightarrow G_{m,2m}$$
 $G_{m,2m}$ being the Grassmannian variety of $m$-dimensional sub
 spaces of $K^{2m}$.

\vs.2cm\ni\textbf{Fact 2.} Indexing the simple roots of $G$ as in
\cite{bour}, we
 have an identification $$G/P\cong L_m$$ $P$ being the maximal
 parabolic sub group of $G$ corresponding to ``omitting" the
 simple root $\alpha_m$ (the right-end root in the Dynkin diagram of
 $G$).

 \vs.2cm\ni\textbf{Fact 3.} The opposite big
cell $O^{-}$ in $G_{m,2m}$ may be identified as $$O^-=\left\{
\begin{pmatrix}
I_{m} \\ Y
\end{pmatrix},\,Y\in M_{m,m}(K)\right\}\leqno{(*)}$$ where $I_{m}$
is the identity $m\,\times\,m$ matrix, and $M_{m,m}$ is the space
of $m\,\times\,m$ matrices (with entries in $K$). For our purpose,
it will be more convenient to replace $Y$ by $JY$ in the above
identification, where $J$ is the $m\,\times\,m$ matrix with $1$'s
on the anti-diagonal and $0$'s elsewhere.  Then (*) gives rise to
an identification of the opposite big cell $O^{-}_G$ in $L_{m}$ as
$$O^-_G=\left\{
\begin{pmatrix}
I_{m} \\ X
\end{pmatrix},\,X\in Sym\,M_m\right\}\leqno{(**)}$$ where $Sym\,M_m$ is
the space of symmetric $m\,\times\,m$ matrices (with entries in
$K$).

See \cite{g/p-2}) for details.

\subsection{The functions $f_{\tau,\varphi}$ on $O^-_G$ :}\label{opp}
Let $\uj\in I(m,2m), p_{\uj}$ the corresponding Pl\"ucker
co-ordinate on $G_{m,2m}$. Denoting by $f_{\uj}$, the restriction
of $p_{\uj}$ to $O^-$, we have (under the identification (*)), if
$y\in O^-$ corresponds to the matrix $Y$, then
 $f_{\uj} (y)$ is simply the following minor of $Y$: let
 $\uj=(j_1,\dots,j_m)$, and let $j_r$ be the largest
entry $\le m$. Let $\{k_1,\dots,k_{m-r}\}$ be the complement of
$\{j_1,\dots,j_r\}$ in $\{1,\dots,m\}$. Then $f_{\uj}(y)$ is the
$(m-r)$-minor of $Y$ with column indices $k_1,\dots k_{m-r}$, and
row indices $j_{r+1},\dots,j_m$ (here the rows of $Y$ are indexed
as $m+1,\dots,2m$). Conversely, given a minor of $Y$, say,  with
column indices $b_1,\dots,b_s$, and row indices
$j_{m-s+1},\dots,j_m$ (again, the rows of $Y$ are indexed as
$m+1,\dots,2m$), it is $f_{\uj}(y)$, where $\uj=(j_1,\dots,j_m)$
is given as follows: $\{j_1,\dots,j_{m-s}\}$ is the complement of
$\{b_1,\dots,b_s\}$ in $\{1,\dots,m\}$, and $j_{d-s+1},\dots,j_m$
are simply the row indices (see \cite{g/p-2} for details).

\vs.2cm\ni\textbf{The partial order $\ge$:} Given $A=(a_1,\cdots
,a_s),\ A'=(a'_1,\cdots ,a'_{s'}),$ for some $s,s'\ge 1$, we
define
 $A\ge A'$ if
$s\le s',\ a_j\ge a'_j,\ 1\le j\le s$.

We have (cf. \cite{g/p-2}) that on $G/P$, any $p_{\tau,\varphi}$
($(\tau,\varphi)$ being an admissible pair) is the restriction of
some Pl\"ucker co-ordinate on $G_{m,2m}$. Let us denote by
$f_{\tau,\varphi}$ the restriction of $p_{\tau,\varphi}$ to
$O^-_G$. Given $z\in O^-_G$, let $X$ be the corresponding matrix
in Sym$\,M_m$ (under the identification (**)); then as above,
$f_{\tau,\varphi}(z)$  is a certain minor of $X$, say with row
(resp. column) indices $A:=\{a_1,\cdots,a_s\}$ (resp.
$B:=\{b_1,\cdots,b_s\}$); we have, $A\ge B$. Denote this minor by
$p(A,B)$. Conversely, such a minor corresponds to a unique
$f_{\tau,\varphi}$ (see \cite{g/p-2} for details). Thus we have a
bijection $$\theta:\{\mathrm{admissible\ pairs }\}\
{\buildrel{bij}\over{\rightarrow}}\ \{\mathrm{minors\
}p(A,B)\mathrm{\ of\ }X, A\ge B\}$$ $X$ being a symmetric $m\times
m$ matrix of indeterminates. The bijection $\theta$ is described
in more detail in \S \ref{12.9.4} below.

\ni {\bf Convention.} If $\tau=\varphi=(1,\dots,m)$, then
$f_{\tau,\varphi}$ evaluated at $z$ is $1$; we shall make it
correspond to the minor of $X$ with row indices (and column
indices) given by the empty set.

\subsection{The opposite cell in $X(w)$}\label{12.1.10}
For a Schubert variety $X(w)$ in $L_{m}$, let us denote $O^-_G\cap
X(w)$ by $Y(w)$. We consider $Y(w)$ as a closed subvariety of
$O^-_G$. In view of Proposition \ref{present}, we obtain that the
ideal defining $Y(w)$ in $O^-_G$ is generated by
$$\{f_{\tau,\varphi}\mid,\ w\not\ge\tau\}.$$

\subsection{Symmetric Determinantal Varieties}\label{det}

Let $\mathcal{Z}=Sym\,M_{m}$, the space of all symmetric $m \times
m$ matrices with entries in $K$. We shall identify $\mathcal{Z}$
with $\mathbb{A}^{N}$, where $N=\frac{1}{2}m(m+1)$. We have
$K[\mathcal{Z}]=K[z_{i,j},\ 1\le i\le j\le m]$.

\ni{\bf{The variety $D_t(Sym\,M_m)$.}}  Let $X=(x_{ij})$, $1\le i,
 j\le m,x_{ij}=x_{ji}$ be a $m\times m$ symmetric matrix of indeterminates. Let
$A,B, \ A\subset \{1,\cdots ,m\}, \ B \subset \{1,\cdots ,m\}, \
\#A=\#B=s$, where $s\le m$. We shall denote by $p(A,B)$ the
$s$-minor of $X$ with row indices given by $A$, and column indices
given by $B$. For $t,\ 1\le t\le m$, let $I_t$ be the ideal in
$K[\mathcal{Z}]$ generated by $\{p(A,B),\ \ A\subset \{1,\cdots
,m\}, \ B\subset \{1,\cdots ,m\}, \ \#A=\#B=t\}$. Let
$D_t(Sym\,M_m)$
 be the {\em{ symmetric determinantal variety
}}(a closed subvariety of $\mathcal{Z}$), with $I_t$ as the
defining ideal. In the discussion below, we also allow $t=m+1$ in
which case $D_t(Sym\,M_m)=\mathcal{Z}$.

\ni{\bf{Identification of $D_t(Sym\,M_m)$ with $Y(\phi)$.}} Let
$G=Sp_{2m}(K)$.  As in \S \ref{12.1.6}, let us identify the
opposite cell $O^-_G$ in $G/P (\cong L_m)$ as
$$O^-_G=\left\{
\begin{pmatrix}
I_{m} \\ X
\end{pmatrix}\right\}$$ where $X$ is a symmetric $m\times m$
matrix. As seen above (cf. \S \ref{opp}), we have a bijection:
 $$\{f_{\tau,\varphi},\, (\tau,\varphi) \in {\mathcal{A}}\}
 {\buildrel{\texttt{bij}}\over{\rightarrow}}
 \{\texttt{ minors}\, p(A,B), A,B\in I(r,m),A\ge B,0\le r\le m
\texttt{ of }X\}$$ (here, ${\mathcal{A}}$ is the set of all
admissible pairs, and $I(r,m)=\{\ui=(i_1,\dots,i_r)|1\le
i_1<\dots<i_r\le m\}$; also note that as seen in \S \ref{opp}, if
$\tau=\varphi = (1,2, \cdots ,m)$, then $f_{\tau,\varphi}=$ the
constant function $1$ considered as the minor of $X$ with row
indices (and column indices) given by the empty set).


 Let $\phi$ be the $m$-tuple, $\phi = (t,t+1,
\cdots , m, 2m+2-t,2m+3-t, \cdots ,2m)(=(t,t+1, \cdots ,m,
(t-1)',(t-2)',\cdots,1'))$ (note that $\phi$ consists of the two
blocks $[t,m]$, $[2m+2-t,2m]$ of consecutive integers - here, for
$i<j$, $[i,j]$ denotes the set

\ni $\{i,i+1, \cdots , j\}$, and for $1\le s\le 2m,s'=2m+1-s$). If
$t=m+1$, then we set $\phi = (m',(m-1)', \cdots ,1')$; note then
that $Y(\phi) = O^-_G(\cong Sym\,{M_m})$.

\begin{thm}\label{iso}(cf.\cite{g/p-2})
$D_t(Sym\,M_m)\cong Y(\phi) $; further,
dim$\,D_t(Sym\,M_m)=\frac{1}{2}(t-1)(2m+2-t)=
 \frac{1}{2}(t-1)(t-1)'$.
\end{thm}

\begin{cor}\label{12.9.1.8}
$K[D_t(Sym\,M_m)]\cong R(\phi)_{(p_{id})}$, the homogeneous
localization of $R(\phi)$ at $p_{id}$ ($id$ being the $m$-tuple
$(1,\cdots,m)$).
\end{cor}

\vs.2cm\ni\textbf{Singular locus of $D_t(Sym\,M_m)$:} For our
discussion in \S \ref{normal} (especially, in the proof of Lemma
\ref{regular}), we will be required to know Sing$\,D_t(Sym\,M_m)$,
the singular locus of $D_t(Sym\,M_m)$. Let $\phi=([t,m],
[2m+2-t,2m])$ as above. From \cite{g/p-7}, we have
$$Sing\,X(\phi)=X(\phi')$$ where $\phi'=([t-1,m],
[2m+3-t,2m])$. This together with Theorem \ref{iso} implies the
following theorem:
\begin{thm}\label{sing} Sing$\,D_t(Sym\,M_m)=D_{t-1}(Sym\,M_m)$
\end{thm}

\subsection{The set $H_m$}\label{12.9.4} Let
$$H_m=\{(A,B)\in{\underset{0\leq s\leq m}{\cup }}I(s,m)\times
I(s,m)\,|\,A\ge B\}$$ where our convention is that
$(\emptyset,\emptyset)$ is the element of $H_m$ corresponding to
$s=0$. We define $ \succeq$ on $H_m$ as follows:

$\bullet$\hskip.5cm We declare $(\emptyset,\emptyset)$ as the
largest element of $H_m$.

$\bullet$\hskip.5cm For $(A,B),(A',B')$ in $H_m$, say,
$A=(a_1,\cdots ,a_s),\ B=(b_1,\cdots ,b_s),\ A'=(a'_1,\cdots
,a'_{s'}),\ B'=(b'_1,\cdots ,b'_{s'})$ for some $s,s'\ge 1$, we
define
 $(A,B) \succeq (A',B')$ (or also $p(A,B) \succeq p(A',B')$) if
$B\ge A'$ (here, $\ge$ is as in \S \ref{opp})


\ni\textbf{The bijection $\theta$:} Let $X=(x_{ij})$ be a generic
symmetric $m\times m$ matrix. Let $(\tau,\varphi)\in
{\mathcal{A}}$, i.e., $(\tau,\varphi)$ is an admissible pair. As
elements of $I_G(m,2m)$, let
$$\tau=(a_1,\cdots ,a_r,b'_1,\cdots ,b'_{s}),\varphi=(c_1,\cdots
,c_r,d'_1,\cdots ,d'_{s})$$ where $r+s=m,a_i,b_j,c_i,d_j \mathrm{\
are\ }\le n$, and (recall that) for $1\le q\le 2m, q'=2m+1-q$. We
would like to remark that the fact that $(\tau,\varphi)$ is an
admissible pair implies that number of entries $\le m$ in $\tau$
and $\varphi$ are the same (see \cite{g/p-2,g/p-7} for details).
Denote
$$\ui :=(i_1,\cdots,i_m):=(a_1,\cdots ,a_r,d'_1,\cdots
,d'_{s})\,(\in I(m,2m))$$ (here, $I(m,2m)$ is as in \S \ref{jm}).
 Set
$$A_{\ui}=\{2m+1-i_m,2m+1-i_{m-1},\cdots ,2m+1-i_{r+1}\},$$
$$B_{\ui}=\text{ the complement of }\{i_1,i_2\cdots ,i_r\} \text{
in } \{1,2,\cdots ,m\}.$$ Define $\te:{\mathcal{A}} \to
\{\text{all minors of }X\}$ by setting $\te
(\ui)=p(A_{\ui},B_{\ui})$ (here, the constant function $1$ is
considered as the minor of $X$ with row indices (and column
indices) given by the empty set). Then $\te $ is a bijection (cf.
\cite{g/p-2}). Note that $\te$ reverses the respective (partial)
orders, i.e., given $\ui,\ui ' \in I(m,2m)$, corresponding to
admissible pairs $(\tau,\varphi),(\tau',\varphi')$  we have, $\ui
\le \ui ' \iff \te (\ui) \succeq \te (\ui ')$. Using the
comparison order $\succeq$, we define {\em{standard monomials}} in
$p(A,B)$'s for $(A,B)\in H_m$:

\begin{defn}
A monomial $p(A_1,B_1)\cdots p(A_s,B_s), s\in\mathbb{N}$ is said to
be standard
if $p(A_1,B_1)\succeq\cdots \succeq p(A_s,B_s)$.
\end{defn}
In view of  Theorems \ref{basis}, \ref{iso},  we obtain

\begin{thm}\label{iden} Let $H_{t-1}=\{(A,B)\in H_m\,|\,\#A\le t-1\}$.
Standard monomials in $\{p(A,B),\,(A,B)\in H_{t-1}\}$ form a basis
for $K[D_t(Sym\,M_m)]$, the algebra of regular functions on
$D_t(Sym\,M_m)$ ($\subset Sym\,M_m$).
\end{thm}

As a direct consequence of Proposition \ref{qualitative}, we
obtain
\begin{prop}\label{qualitative2} Let $p(A_1,A_2),
p(B_1,B_2)$ (in $K[D_t(Sym\,M_m)]$) be non standard. Let
$$p(A_1,A_2)p(B_1,B_2)=\sum\, a_ip(C_{i1},C_{i2})p(D_{i1},D_{i2}), a_i\in K^*\eqno{(*)}$$ be
  the straightening relation, i.e., R.H.S. is a sum of standard monomials. Then for every $i$,
  $C_{i1},C_{i2},D_{i1},D_{i2}$ have cardinalities $\le t-1$; further,
\begin{enumerate}
\item $C_{i1}\ge$ both $A_1$ and $B_1$; further, if for some $i$,
$C_{i1}$ equals $A_1$ (resp. $B_1$), then $C_{i2}>A_2$ (resp.
$>B_2$) . \item  $D_{i2}\le$ both $A_2$ and $B_2$; further, if for
some $i$, $D_{i2}$ equals $A_2$ (resp. $B_2$), then $D_{i1}<A_1$
(resp. $<B_1$).  \item Suppose there exists a $\sigma\in S_4$ (the
symmetric group on $4$ letters) such that
  $\sigma(A_1)\ge\sigma(A_2)\ge\sigma(B_1)\ge\sigma(B_2)
  $, then $(\sigma(A_1),\sigma(A_2)),(\sigma(B_1),\sigma(B_2))
  $ are in $H_{t-1}$, and $p(\sigma(A_1),\sigma(A_2))p(\sigma(B_1),\sigma(B_2))$
occurs with coefficient $\pm 1$ in (*).
\end{enumerate}

\end{prop}

\begin{rem} On the R.H.S. of (*), $C_{i1},C_{i2}$ could both be the empty set (in
which case $p(C_{i1},C_{i2})$ is understood as $1$). For example,
with $X$ being a $2\times 2$ symmetric matrix of indeterminates,
we have
$$p_{2,1}^2=p_{2,2}p_{1,1}-p_{\emptyset,\emptyset}p_{12,12}.$$
\end{rem}
\begin{rem}\label{homo}
In the sequel, while writing a straightening relation as in
Proposition \ref{qualitative2}, if for some $i$, $C_{i1},C_{i2}$
are both the empty set, we keep the corresponding
$p(C_{i1},C_{i2})$ on the right hand side of the straightening
relation (even though its value is $1$) in order to have
homogeneity in the relation.
\end{rem}

 Taking $t=m+1$ (in which case $D_t(Sym\,M_m)=Z=Sym\,M_{m}$) in Theorem \ref{iden}
and Proposition \ref{qualitative2}, we obtain

\begin{thm}\label{iden1}
\begin{enumerate}
\item Standard monomials in $\{p(A,B),\,(A,B)\in H_m\}$'s form a
basis for $K[\mathcal{Z}](\cong K[x_{ij},1\le i\le j\le m] )$ (if
$(A,B)=(\emptyset,\emptyset)$, then $p(A,B)$ should be understood
as the constant function $1$). \item Relations similar to those in
Proposition \ref{qualitative2} hold on $\mathcal{Z}$.
\end{enumerate}
\end{thm}
\begin{rem} Note that Theorem \ref{iden} recovers Theorem 5.1 of
\cite{d-p}. But we had taken the above approach of deducing
Theorem\ref{iden}  from Theorems \ref{basis}, \ref{iso} in order
to derive the straightening relations as given by Proposition
\ref{qualitative2} (which are crucial for the discussion in \S
\ref{stan}).
\end{rem}

\vs.2cm\ni{\bf{A presentation for $K[D_t(Sym\,M_m)]$.}} Consider
the polynomial algebra

\ni $K[x(A,B), (A,B)\in H_{t-1}]$. For two non-comparable pairs
(under $\succ$ (cf. \S\ref{12.9.4}))

\ni $(A_1,A_2), (B_1,B_2)$ in $H_{t-1}$, denote

$$F((A_1,A_2);(B_1,B_2))=x(A_1,A_2)(B_1,B_2)-\sum\,
a_ix(C_{i1},C_{i2})x(D_{i1},D_{i2})$$ where

\ni $C_{i1},C_{i2}, D_{i1},D_{i2},a_i$ are as in Proposition
\ref{qualitative2}. Let $J_{t-1}$ be the ideal generated by
$$\{F((A_1,A_2);(B_1,B_2)), (A_1,A_2), (B_1,B_2) \texttt{ non-comparable}\}$$ Consider the surjective map $f_{t-1}:K[x(A,B), (A,B)\in
H_{t-1}]\rightarrow K[D_t(Sym\,M_m)]$, $x(A,B)\mapsto p(A,B)$.
Then in view of Proposition \ref{present} and Theorem \ref{iso},
we obtain

\begin{prop}\label{present2}
$ f_{t-1}$ induces an isomorphism
$$K[x(A,B), (A,B)\in H_{t-1}]/J_{t-1}\cong
K[D_t(Sym\,M_m)]$$
\end{prop}

\section{$O_{n}(K)$ actions}\label{son}
\setcounter{subsubsection}{0} In this section, we give a
GIT-theoretic proof of the \textbf{First and Second fundamental
theorems} for $O_n(K)$-actions appearing in classical invariant
theory (cf.\cite{d-p,weyl}). Let $V=K^{n}$, together with a
non-degenerate,
 symmetric bilinear form $\langle,\rangle$. Let $G=O(V)$ (the orthogonal group consisting of linear automorphisms of
$V$ preserving $\langle,\rangle$). We shall take the matrix of the
form $\langle,\rangle$ to be
$J_n:=\texttt{anti-diagonal}(1,\cdots,1)$. Then a matrix $A\in
GL(V)$ is in $O(V)$ if and only if
$$^tAJ_nA=J_n,\ i.e,\ J_n^{-1}(^tA^{-1})J_n=A,\ i.e,\ J_n(^tA^{-1})J_n=A\
(\texttt{note that }J_n^{-1}=J_n)$$

\begin{rem}\label{char} In particular, note that a diagonal matrix
$A=\texttt{diagonal}(t_1,\cdots,t_n)$ is in $O(V)$ if and only if
$t_{n+1-i}=t_i^{-1}$.
\end{rem}
 Let  $X$ denote the affine space
$V^{\oplus m}=V\oplus \cdots \oplus V,$ where $m>n.$ For
${\underline{u}}=(u_{1},u_{2},...,u_{m})\in X$, writing
$u_i=(u_{i1},\cdots,u_{in})$ (with respect to the standard basis
for $K^n$), we shall identify ${\underline{u}}$ with the $m\times
n$ matrix $U:=(u_{ij})_{m\times n}$. Thus, we identify $X$ with $
M_{m,n},$ the space of $m\times n$ matrices (with entries in $K$).
Consider the diagonal action of $O(V)$ on $X$. The induced action
on $K[X]$ is given by

\begin{equation*}
(A\cdot f)({\underline{u}})=f(A\cdot{\underline{u}}
)=f(UA),\,\ui\in X,\,f\in K[X],\,A\in O_n(K)(=O(V))
\end{equation*}

\vs.2cm\ni\textbf{A basic Lemma on quotients:} Consider a linear
action of a reductive group $G$ on an affine variety $X=SpecR$
with $R$ a graded $K$-algebra. Let $f_1,\cdots,f_N$ be homogeneous
$G$-invariant elements in $R$. Let $S=K[f_1,\cdots,f_N]$. We
recall (cf. \cite{g/p-2,l-s}) the Lemma below which gives a set of
sufficient conditions for the equality $S=R^G$.

 Let $X^{ss}$ be the set
of semi-stable points of $X$ (i.e., points $x$ such that $0\not\in
{\overline{G\cdot x}}$). Let $\psi:X\rightarrow\mathbb{A}^N$ be the
map, $x\mapsto(f_1(x),\cdots,f_N(x))$. Denote $D=SpecS$. Then $D$
is the categorical quotient of $X$ by $G$ and $\psi:X\to D$ is the
canonical quotient map, provided the following conditions are
satisfied:
\begin{lem}\label{normality} (cf. \cite{g/p-2,l-s})
(i) For $x\in X^{\text{\rm ss}}$, $\psi(x)\ne (0)$.

(ii) There is a $G$-stable open subset $U$ of $X$ such that $G$
operates freely on $U,\,U\rightarrow U/G$ is a $G$-principal fiber
space, and $\psi$ induces an immersion $U/G\to\mathbb{A}^N$ (i.e. an
injective morphism with the induced maps between tangent spaces
injective).

(iii) $\dim D=\dim U/ G$.

(iv) $D$ is normal.

Then $\psi:X\to D$ is the categorical quotient of $X$ by $G$.
\end{lem}

\ni\textbf{The functions $\varphi _{ij}$:} Consider the functions
$\varphi _{ij}:X\rightarrow K$ defined by $ \varphi
_{ij}(({\underline{u}} ))=\langle u_i,u_j\rangle$, $1\leq i,j\leq
m$. Each $\varphi _{ij}$ is clearly in $K[X]^{O(V)}$.

Let $S$ be the subalgebra of $K[X]^{O(V)}$ generated by $\{\varphi
_{ij}\}$. We shall now show (using Lemma \ref{normality}) that $S$
is in fact equal to $K[X]^{O(V)}$.

\begin{thm}\label{fund}
The morphism $\psi:X\rightarrow Sym \,M_m,
\,{\underline{u}}\mapsto \left(\langle u_i,u_j\rangle\right)$ is
$O(V)$-invariant. Furthermore, it maps $X$ onto $D_{n+1}(Sym\,M_m)$
and identifies the categorical quotient $X\fs O(V)$ with
$D_{n+1}(Sym\,M_m)$ (here, $D_{n+1}(Sym\,M_m)$ is as in \S
\ref{det} with $t=n+1$).
\end{thm}
\begin{rem}\label{roandso}
In some parts of the proof below of the above theorem,  it is
assumed that $m\geq n$.    This is not however a serious issue.
Suppose now that $m<n$.  Then the proof can be adapted to this case.
Furthermore, the proof with obvious modifications
shows also that the theorem holds
with $O(V)$ replaced by $SO(V)$.    Observe that 
the variety $D_{n+1}(Sym\,M_m)$ is the same as the affine space $Sym\,M_m$.
Thus $X\fs O(V)=X\fs SO(V)\cong Sym\,M_m$.  In other words,
$K[X]^{O(V)}=K[X]^{SO(V)}$ is a polynomial algebra.
\end{rem}
\begin{proof}
Clearly, $\psi(X) \subseteq D_{n+1}(Sym\,M_m)$ (since,
$\psi(X)=Spec\,S$, and clearly $Spec\,S\subseteq
D_{n+1}(Sym\,M_m)$ (since any $n+1$ vectors in $V$ are linearly
dependent)). We shall prove the result using Lemma
\ref{normality}. To be very precise, we shall first check the
conditions (i)-(iii) of Lemma \ref{normality} for
$\psi:X\rightarrow Sym\,M_{m}$, deduce that the inclusion
$Spec\,S\subseteq D_{n+1}(Sym\,M_m)$ is in fact an equality, and
hence conclude the normality of $Spec\,S$ (condition (iv) of Lemma
\ref{normality}).

(i) Let $x={\underline{u}}=(u_1,\dots,u_m)\in X^{\text{\rm ss}}$.
Let $W_x$ be the subspace of $V$ spanned by $u_i$'s. Let
$r=dim\,W_x$. Then $r>0$ (since $x\in X^{\text{\rm ss}}$). Assume
if possible that $\psi(x)=0$, i.e. $\langle u_i,u_j\rangle=0$ for
all $i,j$. This implies in particular that $W_x$ is totally
isotropic; hence, $r\le [\frac{n}{2}]$, integral part of
$\frac{n}{2}$. Hence we can choose a basis $\{e_1,\cdots ,e_n\}$
of $V$ such that $W_x=$ the $K$-span of $\{e_1,\cdots ,e_r\}$.
Writing each vector $u_i$ as a row vector (with respect to this
basis), we may represent ${\underline{u}}$  by the $m\times n$
matrix $\mathcal{U}$ given by
$$\mathcal{U}:=\begin{pmatrix}
u_{11}&u_{12}&\dots&u_{1r}&0&\dots&0\\
u_{21}&u_{22}&\dots&u_{2r}&0&\dots&0\\
\vdots&\vdots&&\vdots&\vdots&\vdots&\vdots\\
u_{m1}&u_{m2}&\dots&u_{mr}&0&\dots&0
\end{pmatrix}$$
Choose integers $a_1,\cdots ,a_r, a_{r+1},\cdots, a_n$, so that
$a_i>0,i\le [\frac{n}{2}]$, and $a_{n+1-i}=-a_i,i\le
[\frac{n}{2}]$ (if $n$ is odd, say, $n=2\ell +1$, then, we take
$a_{\ell +1}$ to be $0$).

\vs.1cm\ni Let $g_t$ be the diagonal matrix $ g_t=
diag(t^{a_1},\cdots t^{a_r},t^{a_{r+1}},\cdots , t^{a_{n}})$ (note
that $g_t\in O(V)$ (cf. Remark \ref{char})). Consider the one
parametric subgroup $\{g_t,t\in K^*\}$. We have,
$g_tx=g_t\cdot\mathcal{U}=\mathcal{U}g_t$ $=\mathcal{U}_t$, where
$$\mathcal{U}_t=\begin{pmatrix}t^{a_1}u_{11}&t^{a_2}u_{12}&\dots&t^{a_r}u_{1r}&0&\dots&0\\
t^{a_1}u_{21}&t^{a_2}u_{22}&\dots&t^{a_r}u_{2r}&0&\dots&0\\
\vdots&\vdots&&\vdots&\vdots&\vdots&\vdots\\
t^{a_1}u_{m1}&t^{a_2}u_{m2}&\dots&t^{a_r}u_{mr}&0&\dots&0\end{pmatrix}$$
Hence $g_tx\rightarrow 0$ as $t\rightarrow 0$ (note that $r\le
[\frac{n}{2}]$, and hence $a_i>0, i\le r$), and this implies that
$0\in {\overline{G\cdot x}}$ ($G$ being $O_n(K)$) which is a
contradiction to the hypothesis that $x$ is semi-stable. Therefore
our assumption that $\psi(x)=0$ is wrong and (i) of Lemma
\ref{normality} is satisfied.

(ii)  Let $$U=\{{\underline{u}}\in X\mid
\{u_1,\dots,u_n\}\text{\rm are linearly independent}\}$$ Clearly,
$U$ is a $G$-stable open subset of $X$.

\ni\textbf{Claim :} $G$ operates freely on $U$, $U\rightarrow
U\,mod\,G$ is  a $G$-principal fiber space, and $\psi$ induces an
immersion $U/G\to Sym\,M_{m}$.

\ni\textbf{Proof of Claim:} Let $H=GL_n(K)$. We have a
$G(=O_n(K))$-equivariant identification $$U\cong
H\times\underbrace{V\times\dots\times V}_{(m-n)\,\text{\rm
copies}}=H\times F,\ \mathrm{say}\leqno{(*)}$$ where
$F=\underbrace{V\times\dots\times V}_{(m-n)\,\text{\rm copies}}$.
From this it is clear that $G$ operates freely on $U$. Further, we
see that $U\,mod\,G$ may be identified with the fiber space with
base $H\, mod\, G$, and fiber $\underbrace{V\times\dots\times
V}_{(m-n)\,\text{\rm copies}}$ associated to the principal fiber
space

\ni $H\rightarrow H\,/\,G$. It remains to show that $\psi$ induces
an immersion $U/G\to\mathbb{A}^N$, i.e., to show that the map
$\psi:U/G\to \mathbb{A}^N$ and its differential $d\psi$ are both
injective. We first prove the injectivity of
$\psi:U/G\to\mathbb{A}^N$. Let $x,x'$ in $U/G$ be such that
$\psi(x)=\psi(x')$. Let $\eta,\eta'\in U$ be lifts for $x,x'$
respectively. Using the identification (*) above, we may write
$$\begin{gathered}
\eta=(A,u_{n+1},\cdots,u_m),\,A\in
H\\
\eta'=(A',u'_{n+1},\cdots,u'_m),\,A'\in H
\end{gathered}$$ (here, $u_i,1\le i\le n$, are given by the rows of $A$,
while $u'_i,1\le i\le n$, are given by the rows of $A'$. The
hypothesis that $\psi(x)=\psi(x')$ implies in particular that
$$\langle u_i,u_j\rangle=\langle u'_i,u'_j\rangle\ , 1\le i,j\le
n\ $$ which may be written as $$AJ_n\,\,^tA=A'J_n\,\,^tA'$$ where
$J_n$ is the matrix of the form $\langle ,\rangle$ (note that
since we are writing a vector $v\in V$ as a row vector, $\langle
v_i,v_j\rangle=v_iJ_n\,\,^tv_j$). Hence we obtain
$$(A'^{-1}A)J_n\,\,^t(A'^{-1}A)=J_n,\  i.e.,\  A'^{-1}A\in G$$ This implies that
$$A=A'\cdot g,\ \mathrm{for\ some\ g\in G}\leqno{(**)}.$$    Hence on $U/G$,
we may suppose that
$$\begin{gathered}x=(u_1,\cdots,u_n,u_{n+1},\cdots,u_m)\\
x'=( u_1,\cdots,u_n,u'_{n+1},\cdots,u'_m)
\end{gathered}$$ where $\{u_1,\cdots,u_n\}$ is linearly
independent.

 For a given $j, n+1\le j\le
m$, we have,
$$\langle u_i,u_j\rangle=\langle  u_i,u_j'\rangle\ , 1\le i\le
n,\
 \text{\rm  implies}, u_j=u_j'$$ (since, $\{u_1,\cdots,u_n\}$ is linearly
independent and the form $\langle ,\rangle$ is non-degenerate).
  Thus we obtain $$u_j=u_j',\,{\text{for\  all\  }}j\leqno{(\dagger)}.$$

  The injectivity of $\psi:U/G\to\mathbb{ A}^N$ follows from
  ($\dagger$).

  To prove that
the differential d$\psi$ is injective, we merely note that the
above argument remains valid for the points over $K[\e]$, the
algebra of dual numbers ($=K\oplus K\e$, the $K$-algebra with one
generator $\e$, and one relation $\e^2=0$), i.e., it remains valid
if we replace $K$ by $K[\e]$, or in fact by any $K$-algebra.

(iii) We have $$\dim U/G=\dim U-\dim G=mn-\frac{1}{2}n(n-1)=\dim
D_{n+1}(Sym\,M_{m})$$ (cf.Theorem \ref{iso} with $t=n+1$). The
immersion $U/G\hookrightarrow Spec\, S (\subseteq
D_{n+1}(Sym\,M_{m}))$ together with the fact above that $\dim
U/G=\dim D_{n+1}(Sym\,M_{m})$ implies that $Spec\, S$ in fact
equals $ D_{n+1}(Sym\,M_{m})$.

 (iv) The normality of $Spec\, S (=D_{n+1}(Sym\,M_{m}))$ follows from Theorem \ref{iso}
 (and the normality of Schubert varieties).
\end{proof}
\begin{thm}\label{fund'}
\begin{enumerate}

\item {\bf First fundamental theorem}

 The ring of invariants $K[X]^{O(V)}$ is generated by
$\varphi_{ij}=\langle u_i,u_j\rangle,\ 1\le i,j\le m$.

\item {\bf Second fundamental theorem}

The ideal of relations among the generators in (1) is generated by
the $(n+1)$-minors of the symmetric $m\times m$-matrix
$(\varphi_{ij})$.
\end{enumerate}
\end{thm}
Further, we have, (in view of  Theorem \ref{iden})
\begin{thm}\label{fund''}(\textbf{A standard monomial basis for $K[X]^{O(V)}$})
The ring of invariants $K[X]^{O(V)}$ has a basis consisting of
standard monomials in the regular functions $p(A,B),\, A,B\in
I(r,m),A\ge B,\,r\le n$.
\end{thm}
\begin{rem}\label{repeat}
In view of the above Theorem, we have that the relations given by
Proposition \ref{qualitative2} hold in $K[X]^{O(V)}$.
\end{rem}

\section{The algebra $S$}\label{alg}
\setcounter{subsubsection}{0} Let $V=K^n$ together with a
non-degenerate bilinear form $\langle ,\rangle$. Let $X=V^{\oplus
m} (=V\oplus \cdots \oplus V\,(m\texttt{ copies}))$, and
$G=SO_n(K)$. Denote $R=K[X]$. Our goal is to prove
Cohen-Macaulayness of $R^G$. We accomplish this by proving the
Cohen-Macaulayness of a certain subalgebra $S$ of  $R^G$, and
showing $S$ in fact equals $R^G$. Also, the case $m<n$ is trivial
because the invariant ring in this case is a polynomial ring (see
Remark~\ref{roandso}); for
the case $m=n$, the $K[X]^{O(n)}$  is a polynomial ring, and the
invariant ring $K[X]^{SO(n)}$ is generated by one more element (of
degree two) over $K[X]^{O(n)}$ and is easily seen to be
Cohen-Macaulay. Thus throughout we will suppose $m>n$.

\ni\textbf{The functions $p(A,B)$:} For $A, B $ in $I(r,m),$ where
$ 1\leq r\leq n$, and $A\ge B$, let $p(A,B)$ denote the the
regular function on $X, \,p(A,B)(({\underline{u}} ))=$ the
$r$-minor of the symmetric $m\times m$ matrix $\left( \langle
u_{i},u_{j}\rangle \right) $ with row indices given by the entries
of $A$,  and column indices given by the entries of $B$.

 \ni\textbf{The
functions $u(I)$:} For $I\in I(n,m)$, let $u(I)$ be the function
$u(I):X\rightarrow K, u(I)({\underline{u}}):=$ the $n$-minor of
$U$ with the row indices given by the entries of $I$ (here,
$U=(u_{ij})$ is as in \S \ref{son}). We have, $g\cdot
u(I)=(det\,g)u(I),\,g\in O(V)$. Hence $u(I)$ is in $K[X]^{SO(V)}$.

\begin{lem}\label{formula} For $I,J\in I(n,m)$, we have
$u(I)u(J)=p(I,J)$.
\end{lem}
\begin{proof} Let $M,N$ denote the $n\times n$ submatrices of the $m\times n$
matrix $U=(u_{ij})$ with row indices given by $I,J$ respectively.
We have
$$
u(I)u(J)=(det\,M)\,(det\,N)=(\det\,M)\,(det\,^tN)=det(\,M\,^tN)=p(I,J)
$$ (note that $M\,^tN$ is the submatrix of $\left( \langle
u_{i},u_{j}\rangle \right) $ with row indices given by the entries
of $I$,  and column indices given by the entries of $J$.)
\end{proof}
\vs.2cm\ni\textbf{The algebra $S$:} Let $S$ be the subalgebra of
$R^G$ generated by $p(A,B), A,B \in I(r,m), r\le n-1, A\ge B,\
u(I), I\in I(n,m)$.

\vs.2cm
\begin{rem}
The $K$-algebra $S$ could have been simply defined as the
$K$-subalgebra of $R^G$ generated by $\{\left\langle u_{i},u_j
\right\rangle\}$ (i.e., by $\{p(A,B),\#A=\#B=1\}$) $\cup\
\{u(I),I\in I(n,m)\}$. But we have a purpose in defining it as
above, namely, the standard monomials (in $S$) will be built out
of the $p(A,B)$'s with $\#A\le n-1$, and the $u(I)$'s (cf.
Definition \ref{std} below).
\end{rem}

\subsection{Standard monomials and their linear independence:}\label{stan} In this
subsection, we define standard monomial in $S$, and prove their
linear independence.

\ni\textbf{The set $H$:} Let $$H_{u}:=I(n,m)$$ (note that $H_u$
indexes the $u(I)$'s). Let $$H_{p}=\{(A,B )\in I(r,m)\times
I(r,m),1\leq r\le n-1,A \ge B\}$$ (note that $H_p$ indexes the
$p(A,B)$'s).

 For $(A,B),(C,D)\in H_p$, define
$(A,B)\succeq (C,D)$ as in \S \ref{12.9.4}. Note that ``$\succeq
$" is not a partial order (since $(A,B)\not\succeq (A,B)$, if $A>
B$); nevertheless, it is transitive. It is only a comparison order
(cf. \S \ref{12.9.4} ). Let
$$H=H_p\cup H_u$$ We define a comparison order $\ge$ on $H$ as follows:

(i) For elements $H_p$, it is just the comparison
 order on $H_p$.

 (ii) For elements $H_u$, it is the partial order $\ge$ as defined
 in \S \ref{opp}.

(iii) No element of $H_u$ is greater than any element of $H_p$.

(iv) For $(A,B)\in H_p,\,I\in H_u$, $(A,B)\ge I$, if $B\ge I$
(again, $\ge$ being as in \S \ref{opp}).

 \vs.2cm\ni Thus $S$ has a set of algebra
 generators $\{p(A,B), (A,B)\in H_p,u(I), I\in H_u\}$ indexed by
 the comparison-ordered set $H$.

 \begin{defn}\label{std} A monomial
 $$F=p(A_1,B_1)\cdots
p(A_r,B_r)u(I_1)\cdots u(I_s)$$ in $\{p(A,B),u(I), (A,B)\in H_p,
I\in H_u\}$ is said to be \emph{standard} if
$$(A_1,B_1)\ge \cdots \ge (A_r,B_r)\ge I_1\ge \cdots \ge I_s$$
i.e, $A_1\ge B_1\ge A_2\ge\cdots \ge B_r\ge I_1\ge\cdots\ge I_s$.

\end{defn}

\subsection{Linear independence of standard monomials:}\label{ind}
In this subsection, we prove the linear independence of standard
monomials.
\begin{lem}\label{special} Let $(A,B)\in H_p, I\in H_u$.
\begin{enumerate}
\item The set of standard monomials in the $p(A,B)$'s  is linearly
independent. \item The set of standard monomials in the $u(I)$'s
 is linearly independent.
\end{enumerate}
\end{lem}

\begin{proof}
As in the proof of Proposition \ref{relations}, the subalgebra
generated by

\ni $\{p(A,B), \, A,B \in H_n\}$ being  $R^{O(V)}$, gets
identified with $K[D_{n+1}(Sym(\,M_m)]$, and hence (1) follows
from Theorem \ref{iden}.

\ni (2) follows from \cite{l-s}, Theorem 1.6.6,(1) applied to
$K[u_{ij},1\le i\le m,1\le j\le n]$.
\end{proof}

\begin{prop}\label{indpt}
Standard monomials (in $S$)are linearly independent.
\end{prop}
\begin{proof}  Let
$$F=G+H=0,\leqno{(*)}$$ be
a relation among standard monomials, where
$G=\sum\,c_iG_i,\,H=\sum\,d_jH_j$ where for each $p(A_1,B_1)\cdots
p(A_r,B_r)u(I_1)\cdots u(I_s)$ in $G$ (resp. $H$), $s$ is even -
including $0$ - (resp. odd). Consider a $g$  in $O_n(K)$, with
det$\,g=-1$. Then noting that $u(I)u(J)=p(I,J)$, and using the
facts that $g\cdot p(A,B)=p(A,B),g\cdot u(I)=(det\,g) u(I)$, we
have, $F-gF=\sum\,2d_jH_j=0$. Since char$(K)\ne 2$, if we show
that $\sum\,d_jH_j=0$, then it would follow (in view of Theorems
\ref{fund} and \ref{iden}) that (*) is the trivial sum. Thus  we
may suppose
$$F=\sum\,d_jH_j\leqno{(**)}$$ where each $H_j$ is a standard monomial of the form:
$$H_j=p(R_1,S_1)\cdots
p(R_l,S_l)u(I)$$ Now multiplying $(**)$ by $u(I_n),I_n$ being
$(1,\cdots,n)$ (and noting as above the equality $u(I)u(I_n)=p(I,I_n)$),
we obtain
$$\sum d_iP_i=0$$ where each $P_i$ is a standard monomial of the form $p(U_1,Q_1)\cdots
p(U_m,Q_m)$ (note that for each standard monomial $H_j$ appearing
in $(**)$, $H_ju(I_n)$ is again standard). Now the required result
follows from the linear independence of $p(A,B)$'s in
$K[X]^{O(V)}$(cf. Theorems \ref{fund} and \ref{iden}).

\end{proof}

\section{The algebra $S({D})$}\label{alg'} Let $S$ be as in the previous section. To prove the generation
of $S$ (as a $K$-vector space) by standard monomials, we define a
$K$-algebra $S({D})$, construct a standard monomial basis for
$S({D})$ and deduce the results for $S$ (in fact, it will turn out
that $S({D})\cong S$). We first define the $K$-algebra $R({D})$ as
follows:

Let $${D}=H\cup \{\mathbf{1}\}$$ $H$ being as in \S \ref{stan}.
Extend the comparison order on $H$ to ${D}$ by declaring
$\{\mathbf{1}\}$ as the largest element. Let $P({D})$ be the
polynomial algebra
$$P({D}):=K[X(A,B),Y(I),X(\mathbf{1}), (A,B)\in H_p,I\in
H_u ]$$ Let $\mathfrak{a}({D})$ be the homogeneous ideal in the
polynomial algebra $P({D})$ generated by the relations (1)-(3) of
Proposition \ref{relations} ($X(A,B),Y(I)$ replacing $p(A,B),u(I)$
respectively), with relation (2) homogenized as
$$X(A_1,A_2)X(B_1,B_2)=\sum\,
a_iX(C_{i1},C_{i2})X(D_{i1},D_{i2})$$ where $X(C_{i1},C_{i2})$ is
to be understood as $X(\mathbf{1})$ if both $C_{i1},C_{i2}$ equal
the emptyset (cf. Remark \ref{homo}). Let
$$R({D})=P({D})/\mathfrak{a}({D})$$ We shall denote the classes
of $X(A,B),Y(I),X(\mathbf{1})$ in $R({D})$ by

\ni $x(A,B),y(I),x(\mathbf{1})$ respectively.

\vs.2cm\ni\textbf{The algebra $S({D})$:}  Set
$S({D})=R({D})_{(x(\mathbf{1}))}$, the homogeneous localization of
$R({D})$ at $x(\mathbf{1})$. We shall denote
$\frac{x(A,B)}{x(\mathbf{1})}, \frac{y(I)}{x(\mathbf{1})}$ (in
$S({D})$) by $c(A,B), d(I)$ respectively.

Let $\varphi_{D}:S({D})\rightarrow S$ be the map
$\varphi_{D}(c(A,B))= p(A,B), \varphi_{D}(d(I))=u(I)$. Let
$$\theta_{D}:R({D})\rightarrow S({D})$$ be
 the canonical map.
Denote $\gamma_{D}:R({D})\rightarrow S$ as the composite
$\gamma_{D}=\varphi_{D}\circ\theta_{D}$.

\subsection{A standard monomial basis for $R({D})$:} We define a monomial in

\ni $x(A,B),y(I),x(\mathbf{1}))$ (in $R({D})$)to be standard in
exactly the same way as in Definition \ref{std}.
\begin{prop}\label{indpt'} The standard monomials  in the
$x(A,B),y(I),x(\mathbf{1})$ are linearly independent.
\end{prop}

\begin{proof} The result follows by
considering $\gamma_{D}:R({D})\rightarrow S$, and using the linear
independence of standard monomials in $S$ (cf. Proposition
\ref{indpt}).
\end{proof}

\subsection{Quadratic relations:}

Before proving the generation of $R({D})$ by standard monomials,
we first describe certain ``straightening relations'' 
(expressions for non-standard monomials as linear sums of standard
monomials) among $p(A,B)'s,u(I)'s$ (in $S$), to be used while
proving generation of $R({D})$ by standard monomials.
\begin{prop}\label{relations}
\begin{enumerate}
 \item Let $I,I' \in
H_u$ be not comparable. We have,
\begin{equation*}
u(I)u(I')=\underset{r}{\sum }b_{r}u(I_{r})u(I'_{r}),\, b_{r}\in
K^*
\end{equation*}
where for all $r$, $I_r\ge$ both $I,I'$, and $I'_r\le$ both
$I,I'$; in fact, $I_r>$ both $I,I'$ (for, if
$I_r=I\mathrm{\,or\,}I'$, then $I,I'$ would be comparable).

\item Let $(A_1,A_2), (B_1,B_2)\in H_p$ be not comparable. Then we
have
$$p(A_1,A_2)p(B_1,B_2)=\sum\, a_ip(C_{i1},C_{i2})p(D_{i1},D_{i2}), a_i\in K^*,$$
where $(C_{i1},C_{i2}), (D_{i1},D_{i2})$ belong to $H_p$, and
$C_{i2}\ge D_{i1}$; further, for every $i$, we have
\begin{enumerate}
\item $C_{i1}\ge$ both $A_1$ and $B_1$; if $C_{i1}= A_1$ (resp.
$B_1$), then $C_{i2}>A_2$  (resp.$B_2$). \item $D_{i2}\le$ both
$A_2$ and $B_2$; if $D_{i2}= A_2$ (resp. $B_2$), then $D_{i1}<A_1$
(resp.$B_1$).
\end{enumerate}

\item Let $I\in H_u, (A,B)\in H_p$ be such that $B\not\ge I$. We
have,
\begin{equation*}
p(A,B)u(I)=\underset{t}{\sum }%
d_{t}p(A_{t},B_{t})u(I_{t}),\,d_{t}\in K^*
\end{equation*}
where for every $t$, we have, $(A_t,B_t)\in H_p,B_t\ge I_t$.
Further, $A_t\ge$ both $A$ and $I$; if $A_t=A$, then $B_t> B$.
\end{enumerate}
\end{prop}
\begin{proof} In the course of the proof, we will be repeatedly
using the fact that the sub algebra generated by $\{p(A,B), \, A,B
\in H_n\}$ (where recall that $H_n=\{(A,B)\in H_m,\, \#A\le
n\},H_m$ being as in \S ref{12.9.4}) being $R^{O(V)}$ (cf. Theorem
\ref{fund'}), the relations given by Proposition
\ref{qualitative2} hold in $K[X]^{O(V)}$ (cf. Remark
\ref{repeat}). We will also use the basic formula
$u(I)u(J)=P(I,J)$ (cf. Lemma \ref{formula}).

\ni Assertion (1) follows from \cite{l-s}, Theorem 4.1.1,(2).

\ni Assertion (2) follows from Proposition \ref{qualitative2} (cf.
Remark \ref{repeat}).

\ni (3). We have, $p(A,B)u(I)u(I_n)=p(A,B)p(I,I_n)$ ($I_n$ being
$(1,2,\cdots ,n)$). The hypothesis that $B\not\ge I$ implies that
$p(A,B)p(I,I_n)$ is not standard. Hence Proposition
\ref{qualitative2} implies that in $K[D_{n+1}(Sym(\,M_m)]$
$$p(A,B)p(I,I_n)=\sum\, a_ip(C_{i1},C_{i2})p(D_{i1},D_{i2}),
a_i\in K^*$$ where R.H.S. is a standard sum, i.e., $C_{i1}\ge
C_{i2}\ge D_{i1}\ge D_{i2},\forall i$. Further, for every $i$,
$C_{i1}\ge$ both $A$ and $I$; if $C_{i1}=A$, then $C_{i2}>B$;
$C_{i2}\ge$ both $B$ and $I_n$; $D_{i2}\le $ both $B$ and $I_n$
which forces $D_{i2}=I_n$ (note that in view of Proposition
\ref{qualitative2}, all minors in the above relation have size
$\le n$); and hence $\#D_{i1}=n$, for all $i$. Hence
$p(D_{i1},D_{i2})=u(D_{i1})u(I_n)$, for all $i$. Hence canceling
$u(I_n)$, we obtain $$p(A,B)u(I)=\sum\,
a_ip(C_{i1},C_{i2})u(D_{i1}),$$ where $C_{i1}\ge$ both $A$ and
$I$; if $C_{i1}=A$, then $C_{i2}>B$. This proves (3).
\end{proof}


\vs.2cm\ni\textbf{Generation of $R({D})$ by standard monomials:}
We shall now show that any non-standard monomial $F$ in $R({D})$
is a linear sum of standard monomials. Observe that if $M$ is a
standard monomial, then $x(\mathbf{1})^lM$ is again standard;
hence we may suppose $F$ to be:
$$F=x(A_1,B_1)\cdots x(A_r,B_r)y(I_1)\cdots y(I_s)$$

Fix an integer $N$ sufficiently large. To each element $A\in
\cup_{r=1}^n\,I(r,m)$, we associate an $(n+1)$-tuple as follows:
Let $A\in I(r,m)$, for some $r$, $A=(a_1,\cdots ,a_r)$. To $A$, we
associate the $n+1$-tuple ${\overline{A}}=(a_1,\cdots ,a_r,
m,m,\cdots ,m,1)$

  To $F$, we associate the integer
$n_F:$ (and call it the \textit{weight of $F$}) which has the
entries of
${\overline{A_1}},{\overline{B_1}},{\overline{A_2}},{\overline{B_2}},\cdots,
{\overline{A_r}},{\overline{B_r}},{\overline{I_1}},\cdots,
{\overline{I_s}}$ as digits (in the $N$-ary presentation). The
hypothesis that $F$ is non-standard implies that

\ni either $x(A_i,B_i)x(A_{i+1},B_{i+1})$ is non-standard for some
$i\le r-1$, or $x(A_r,B_r)y(I_1)$ is non-standard, or
$u(I_t)u(I_{t+1})$ is non-standard for some $t\le s-1$ is
non-standard. Straightening these using Proposition
\ref{relations}, we obtain that $F=\sum\,a_iF_i$ where
$n_{F_{i}}>n_F$, for all $i$, and the result follows by decreasing
induction on $n_F$ (note that while straightening a degree $2$
relation using Proposition \ref{relations}, if $x(\mathbf{1})$
occurs in a monomial $G$, then the digits in $n_G$ corresponding
to $x(\mathbf{1})$ are taken to be $m,m,\cdots,m$ ($2n+2$-times)).
Also note that the largest $F$ of degree $r$ in $x(A,B)$'s and
degree $s$ in the $y(I)$'s is $x(\{m\},\{m\})^ru(I_0)^s$ (where
$I_0$ is the $n$-tuple $(m+1-n,m+2-n,\cdots ,m)$) which is clearly
standard.

Hence we obtain

\begin{prop}\label{gen'} Standard monomials in
$x(A,B),y(I),x(\mathbf{1}))$
 generate  $R({D})$ as a $K$-vector space.
 \end{prop}

 Combining Propositions \ref{indpt'}, \ref{gen'}, we obtain

 \begin{thm}\label{mainprime}
Standard monomials in $x(A,B),y(I),x(\mathbf{1})$
 give a basis for the $K$-vector space $R({D})$.
\end{thm}
 \subsection{Standard monomial bases for $S({D})$}
Standard monomials in
$c(A,B),d(I)$ in $S({D})$
are defined in exactly the
same way as in Definition \ref{std}.
\begin{thm}\label{main''}
Standard monomials in $c(A,B),d(I)$
 give a basis for the $K$-vector space $S({D})$.
 \end{thm}
 \begin{proof}
The linear independence of standard monomials follows as in the
proof of Proposition \ref{indpt'} by considering
$\varphi_{D}:S({D})\rightarrow S$, and using the linear
independence of standard monomials in $S$ (cf. Proposition
\ref{indpt}).

\ni To see the generation of $S({D})$ by standard monomials,
consider a non-standard  monomial $F$ in $S({D})$, say,
$$F=c(A_1,B_1)\cdots c(A_i,B_i)d(I_1)\cdots d(I_k)$$   Then
$F=\theta_{D}(G)$, where

\ni $G=x(A_1,B_1)\cdots x(A_i,B_i)y(I_1)\cdots y(I_k)$. The
required result follows from Proposition \ref{gen'}.
 \end{proof}

\begin{thm}\label{capMain} Standard monomials in $p(A,B),u(I)$
form a basis for the $K$-vector space $S$.
 \end{thm}

 \begin{proof} We already have
 established the linear independence of standard monomials
 (cf. Proposition \ref{indpt}). The generation by standard
 monomials follows by considering the surjective map
 $\varphi_{D}:S({D})\rightarrow
 S$ and using the generation of $S({D})$ by standard monomials (cf. Theorem
 \ref{main''}).
\end{proof}

 \begin{thm}\label{isom} The map
$\varphi_{D}:S({D})\rightarrow
 S$ is an isomorphism of $K$-algebras.
 \end{thm}
\begin{proof}
Under $\varphi_{D}$, the standard monomials in $S({D})$ are mapped
bijectively onto the standard monomials in $S$. The result follows
from Theorems  \ref{main''} and
  \ref{capMain}.
\end{proof}

\begin{thm}\label{present'}\textbf{A presentation for $S$:}
\begin{enumerate}\item  The
$K$-algebra $S$ is generated by $\{p(A,B),u(I), (A,B)\in H_p,I\in
H_u\}$. \item  The ideal of relations among the generators
$\{p(A,B),u(I)\}$ is generated by  relations (1)-(3) of
Proposition \ref{relations}.
\end{enumerate}
\end{thm}
\begin{proof}
The result follows from Theorem \ref{isom} (and the definition of
$S({D})$).
\end{proof}

\begin{rem}
It will be shown in Theorem \ref{equal} that the inclusion
$S\hookrightarrow R^G$ is in fact an equality.
\end{rem}

\section{Cohen-Macaulayness of $S$}\label{cohen}
\setcounter{subsubsection}{0} In this section, we prove the
Cohen-Macaulayness of $S$ in the following steps:

$\bullet$ \hskip 1cm (i) We prove that $S$ is a doset algebra with
straightening law over a doset $D$ contained in $P\times P$, for a
certain partially ordered set $P$.

$\bullet$ \hskip 1cm (ii) $P$ is a wonderful poset.

$\bullet$ \hskip 1cm (iii) Conclude the Cohen-Macaulayness of $S$
using \cite{d-l}.

Let $P$ be a partially ordered set. Recall
\begin{defn}\label{dos}(cf. \cite{d-l})
A doset of $P$ is a subset $D$ of $P\times P$ such that
\begin{enumerate}
\item  The diagonal  $\Delta (P)\subset D$.
 \item If $(a,b)\in D$, then $a\ge b$.

\item  Let $a\geq b\geq c$ in $P$.

(a) Let $(a,b),(b,c)$ be in $D$. Then $(a,c)\in D$.

(b) Let $(a,c)\in D$. Then $(a,b),(b,c)$ are in $D$.
\end{enumerate}
\end{defn}
\begin{rem}
(i) We shall refer to an element $(\alpha,\beta)$ of $D$ as an
\emph{admissible pair}; $(\alpha,\alpha)$ will be called a
\emph{trivial admissible pair}.

\ni (ii) If $a\geq b\geq c\ge d$ in $P$ are such that $(a,d)$ is
in $D$, then $(b,c)\in D$. This follows from Definition \ref{dos},
3(b).
\end{rem}

\begin{defn}\label{dalg}(cf. \cite{d-l})
A doset algebra with straightening law over the doset $D$ is a
$K$-algebra $E$ with a set of algebra generators $\{x(A,B),
(A,B)\in D\}$ such that
\begin{enumerate}
\item $E$ is graded with $E_0=K$, and $E_r$ equals the $K$-span of
monomials of degree $r$ in $x(A,B)$'s. \item ``Standard monomials"
$x(A_1,B_1)\cdots x(A_r,B_r)$ (i.e., $A_1\ge B_1\ge A_2\ge
B_2\cdots \ge A_r\ge B_r$) is a $K$-basis for $E_r$. \item Given a
non-standard monomial $F:=x(A_1,B_1)\cdots x(A_r,B_r)$, in the
straightening relation $$F=\sum\,a_iF_i\leqno{(*)}$$ expressing
$F$ as a sum of standard monomials, writing

\ni $F_i=x(A_{i1},B_{i1})\cdots x(A_{ir},B_{ir})$ (a standard
monomial) we have the following:

\vs.2cm (a) For every permutation $\sigma$ of
$\{A_{1},B_{1},A_{2},B_{2},\cdots, A_{r},B_{r}\}$, we have,

\ni $\{A_{i1},B_{i1},A_{i2},B_{i2},\cdots, A_{ir},B_{ir}\}$ is
lexicographically greater than or equal to $
\{\sigma(A_1),\sigma(B_1),\sigma(A_2),\sigma(B_2),\cdots,\sigma(A_r),\sigma(B_r)\}$
.

\vs.2cm (b) If there exists a $\tau \in S_{2r}$ such that
$\tau(A_1)\ge\tau(B_1)\ge\tau(A_2)\ge\tau(B_2)\ge\cdots\ge\tau(A_r)\ge\tau(B_r)$,
then the monomial $x(\tau(A_1),\tau(B_1))\cdots
x(\tau(A_r),\tau(B_r))$ occurs on the R.H.S. of (*) with
coefficient $\pm 1$. (We shall refer to this situation as
``$\{A_1, B_1, A_2, B_2,\cdots, A_r, B_r\}$ is totally ordered up
to a reshuffle")

\ni (Note that in 3(b), we have (in view of (3),(4) in Definition
\ref{dos})) that $(\tau(A{_i}),\tau(B{_i}))$'s are admissible
pairs.

\end{enumerate}
\end{defn}

\ni\textbf{The discrete doset algebra $K\{D\}$:} The discrete
doset algebra $K\{D\}$ is the doset algebra with straightening
relations given as follows: Let $F:=x(A_1,B_1)\cdots x(A_r,B_r)$
be a non-standard monomial. Then
$$F=
\begin{cases}x(\tau(A_1),\tau(B_1))\cdots
x(\tau(A_r),\tau(B_r)),\mathrm{\,
 if\,} \tau(A_1)\ge\tau(B_1)\ge\cdots\ge\tau(A_r)\ge\tau(B_r)\\
0,\mathrm{\,otherwise}
\end{cases}
$$
for some $\tau\in S_{2r}$ as above.

\begin{rem}\label{two}
The conditions in (3) are equivalent to the corresponding
conditions for $r=2$.
\end{rem}

Let $K\{P\}$ be \emph{the discrete algebra} over $P$, namely,
\emph{the Stanley-Reisner algebra}, defined as the quotient of the
polynomial algebra $K[x_\alpha, \alpha\in P]$ by the ideal
generated by $\{x_\alpha x_\beta,\,\alpha,\beta\in P $
non-comparable\}.

Recall
\begin{prop}\label{mainprop} (cf. \cite{d-l}, Theorem 3.5)
Let $E$ be a doset algebra with straightening law over a doset $D$
inside $P\times P$.
\begin{enumerate}
\item There exists a sequence $E=B_1,B_2,\cdots, B_r=K\{D\}$ of
doset algebras with straightening law over the doset $D$ such that
there exists a flat family $\mathcal{B}_j, 1\le j\le r-1$ (over
Spec$\,K[t]$) with generic fiber Spec$\,B_j$, and special fiber
Spec$\,B_{j+1}$. \item $E$ is Cohen-Macaulay if $K\{P\}$ is.
\end{enumerate}
\end{prop}

\subsection{A doset algebra structure for $R(D)$:} We first define

\vs.2cm\ni\textbf{The partially ordered set $P$:} Let
$$P:=\cup_{r=1}^n\,I(r,m)\,\cup\,\{\mathbf{1}\}$$ Extend the partial
order on $\cup_{r=1}^n\,I(r,m)$ to $P$ by declaring $\mathbf{1}$
to be the largest element.

Let $D$ be as in \S \ref{alg'}, namely,
$$D=H_p\,\cup\,H_u\,\cup\,\{\mathbf{1}\}$$

\begin{lem}
$D$ is a doset inside $P\times P$.
\end{lem}
\begin{proof}
Clearly, $D\subset P\times P$, and contains $\Delta(P)$. Note that
 $H_u$ (as a subset of $P\times P$) is identified with the
 diagonal $\Delta(H_u)$ (inside
$P\times P$); similarly, $\mathbf{1}$ is identified with
$(\mathbf{1},\mathbf{1})$. Also note  that the non-trivial
admissible pairs in $D$ are among the $(A,B)$'s ($(A,B)\in H_p$).
The remaining conditions in Definition \ref{dos} hold in view of
the results in \S \ref{det} and the results of
\cite{d-l,g/p-4,g/p-5}.
\end{proof}

\begin{prop}\label{dose}
$R(D)$ is a doset algebra with straightening laws over the doset
$D$.
\end{prop}

\begin{proof} Condition (1) in Definition \ref{dalg}
follows from the definition of the $K$-algebra $R(D)$; note that
$R(D)$ has algebra generators $\{x(A,B),(A,B)\in
H_p\},y(I),I\in H_u,x(\mathbf{1})$, indexed by $D$. Note that the
generators $\{y(I),I\in H_u,x(\mathbf{1})\}$ are indexed by
trivial admissible pairs, and the generators indexed by
non-trivial admissible pairs are among $\{x(A,B),(A,B)\in H_p\}$.
Condition (2) in Definition \ref{dalg} follows from Theorem
\ref{mainprime}.

\ni\textbf{Verification of Condition (3):} As in the proof of
Theorem 4.1 of \cite{d-l}, in view of Proposition \ref{relations},
it suffices to verify condition (3) in Definition \ref{dalg} for a
degree two non-standard monomial $F$ (cf. Remark \ref{two}). We
divide the verification into the following cases:

\ni\textbf{Case 1:} Let $F=y(I)y(J)$. In this case, the situation
of 3(b) (in Definition \ref{dalg}) does not exist. The condition
3(a) follows from Proposition \ref{relations},(1).

\ni\textbf{Case 2:} Let $F=x(A,B)y(I)$. In this case again, the
situation of 3(b) (in Definition \ref{dalg}) does not exist (since
$B\not\ge I$, and $I$ can not be $>A\mathrm{\,or\,}B$ - note that
no element of $H_u$ is greater than any element of $H_p$).
Condition 3(a) follows from Proposition \ref{relations},(3).

\ni\textbf{Case 3:} Let $F=x(A_1,B_1)x(A_2,B_2)$. Condition 3(a)
follows from Proposition \ref{relations},(2). Condition 3(b)
follows from Theorem \ref{iso}, and \cite{d-l}, Theorem 4.1
(especially, its proof); here, we should remark that one first
concludes such relations for $Y(\phi)$ ($Y(\phi)$ as in Theorem
\ref{iso}), then for $S$, and hence for $R(D)$ (note that $S$
being $R(D)_{(x(\mathbf{1}))}$ such relations in $S$ imply similar
relations in $R(D)$, since $x(\mathbf{1})^l$ is the largest
monomial in any given degree $l$).
\end{proof}

\begin{cor}\label{cohen'}
$R(D)$ is Cohen-Macaulay.
\end{cor}
\begin{proof}
This follows from Propositions \ref{mainprop}, \ref{dose}. Note
that $P$ is a wonderful poset (in the sense of \cite{d-e-p}), in
fact, a distributive lattice, and hence the discrete algebra
$K\{P\}$ is Cohen-Macaulay (cf. \cite{d-e-p}, Theorem 8.1).
\end{proof}
The above Corollary together with the fact that $S$ is a
homogeneous localization of $R(D)$ implies
\begin{thm}\label{main}
$S$ is Cohen-Macaulay.
\end{thm}

\section{The equality $R^{SO_n(K)}=S$}\label{normal}
\setcounter{subsubsection}{0} We preserve the notation from the
previous sections. In this section, we shall first prove that the
morphism $q:Spec\,R^{SO_n(K)}\longrightarrow Spec\,S$ induced by
the inclusion $S\subseteq R^{SO_n(K)}$ is finite, surjective, and
birational. Then, we shall prove that $Spec\,S$ is normal, and
deduce (using Zariski's Main Theorem) that $q$ is an isomorphism,
thus proving that the inclusion $S\subseteq R^{SO_n(K)}$ is an
equality, i.e., $R^{SO_n(K)}=S$. As a consequence, we will obtain
(in view of Theorem \ref{main}) that $R^{SO_n(K)}$ is
Cohen-Macaulay.

\vs.2cm\ni\textbf{Notation:} In this section, for an integral
domain $A,\kappa(A)$ will denote the quotient field of $A$.

\begin{lem}\label{above} Let $\widetilde{A}$ be a finitely generated $K$-algebra ($K$ being an
algebraically closed field of characteristic different from $2$).
Further, let $\widetilde{A}$ be a domain. Let $\gamma$ be an
involutive $K$-algebra automorphism of $\widetilde{A}$. Let
$A=\widetilde{A}^{\Gamma}$ where $\Gamma\cong
\mathbb{Z}/2\mathbb{Z}$ is the group generated by $\gamma$. We
have,
\begin{enumerate}
\item The canonical map $p:Spec\,\widetilde{A}\longrightarrow
Spec\,A$ induced by the inclusion $A\subset \widetilde{A}$ is a
finite, surjective morphism. \item $\kappa(\widetilde{A})$ is a
quadratic extension of $\kappa(A)$.
\end{enumerate}
\end{lem}
\begin{proof} (1) It is easy to see that a finite set $\mathcal{B}$ of $A$-algebra generators of
$\widetilde{A}$ can be chosen so that they are all eigenvectors of
$\gamma$ corresponding to the eigenvalue $-1$. (This is because
for any $f\in \widetilde{A}$, one has
$f=1/2(f+\gamma(f))+1/2(f-\gamma(f))$.)

Since $a:=f.\gamma(f)\in A$ for any $f\in \widetilde{A}$, each
generator satisfies the equation $x^2+a=0$ (over $A$). Also, since
product of any two generators is an eigenvector corresponding to
the eigenvalue $1$, it follows that $\widetilde{A}$ is generated
as an $A$-module by the set of algebra generators $\mathcal{B}$
together with $1$. Therefore $\widetilde{A}$ is a finite module
over $A$. Hence $p$ is a finite morphism; surjectivity of $p$
follows from the fact (cf. \cite{mum}, ch I.7, Proposition 3) that
a finite morphism $f:Spec\,B\longrightarrow Spec\,A$ of affine
varieties is surjective if and only if $f^*:A \longrightarrow B$
is injective. Assertion (1) follows.

\ni (2) From the discussion in (1), we obtain that the $A$-algebra
$\widetilde{A}$ has algebra generators, say, $\{f_1,\cdots,f_r\}$
such that

$\bullet$ $f_i^2\in A,1\le i\le r$

$\bullet$ $f_if_j\in A,\forall i,j$

Hence we obtain that for every $\alpha\in \widetilde{A}, \alpha^2
\in A$. This implies that every $s\in \kappa(\widetilde{A})$
satisfies a quadratic equation $x^2+a$ over $\kappa(A)$; further,
$\kappa(\widetilde{A})$ is a (finite) separable extension of
$\kappa(A)$ (since char$\,K\not= 2$). Assertion (2) follows from
this.
\end{proof}
\begin{cor}\label{finite}
Let $A=R^{O_n(K)},\widetilde{A}=R^{SO_n(K)}$. We have,
\begin{enumerate}
\item The canonical map $p:Spec\,\widetilde{A}\longrightarrow
Spec\,A$ induced by the inclusion $A\subset \widetilde{A}$ is a
finite, surjective morphism. \item $\kappa(\widetilde{A})$ is a
quadratic extension of $\kappa(A)$.
\end{enumerate}
\end{cor}
\begin{proof}
Taking $\gamma\in O_n(K)$ to be an order $2$ element which
projects onto the generator of
$O_n(K))/SO_n(K)=:\Gamma(=\mathbb{Z}/2\mathbb{Z})$, we have,
$\gamma$ defines an involutive $K$-algebra automorphism of
$\widetilde{A}$, with $\widetilde{A}^{\Gamma}=A$. The result
follows from Lemma  \ref{above}.
\end{proof}
\begin{prop}\label{birational}
The morphism $q:Spec\,R^{SO_n(K)}\longrightarrow Spec\,S$ induced
by the inclusion $S\subseteq R^{SO_n(K)}$ is surjective, finite,
and birational.
\end{prop}
\begin{proof}
Denote $\widetilde{Y}=Spec(R^{SO_n(K)}), Y=Spec(R^{O_n(K)}),
Z=Spec\,S.$ Consider the inclusions $$R^{O_n(K)}\subset S\subseteq
R^{SO_n(K)}\leqno{(*)}$$ Note that the first inclusion is a strict
inclusion, since, $S=R^{O_n(K)}[u(I),I\in I(n,m)]$ (and
$u(I)\not\in R^{O_n(K)}$ for $I\in I(n,m)$). This induces the
following commutative diagram
$$ \commdiag{ \widetilde{Y}&&\cr
\mapdown\lft{q}&\arrow(3,-2)\rt{p}\cr Z&\mapright&Y\cr } $$ The
 finiteness of $q$ follows from the finiteness of $p$
 (cf. Corollary \ref{finite}, (1)); this together with the inclusion
 $q^*:S\hookrightarrow R^{SO_n(K)}$ implies the surjectivity of
 $q$ (cf. \cite{mum}, ch I.7, Proposition 3).

\ni  Now the inclusions given by (*) give the following inclusions
of the respective quotient fields:
$$\kappa(R^{O_n(K)})\subset
\kappa(S)\subseteq \kappa(R^{SO_n(K)})$$ (with the first inclusion
being a strict inclusion). This together with the fact that
$\kappa(R^{SO_n(K)})$ is a quadratic extension of
$\kappa(R^{O_n(K)})$ (cf. Corollary \ref{finite}, (2)) implies
that $\kappa(S)$ is a quadratic extension of $\kappa(R^{O_n(K)})$
as well. Hence we obtain that $\kappa(S)= \kappa(R^{SO_n(K)})$
proving the birationality of $q$.
\end{proof}

 \vs.2cm Finally, to verify the hypotheses in Zariski's Main Theorem for the morphism $q$,
 it remains to show that $S$ is a normal domain. Again in view of
 Theorem \ref{main} and Serre's criterion for normality (namely, $Spec\,A$
is normal if and only if $A$ has $S_2$ and $R_1$), to prove the
normality of $S$, it suffices to show that $Spec\,S$ is regular in
codimension $1$ (i.e., singular locus of $Spec\,S$ has codimension
at least $2$). Towards proving this, we first obtain a criterion
for the branch locus of a finite morphism to have codimension at
least $2$.

  Let $\pi:Z\longrightarrow Y$ be a finite morphism
where $Y$ is a reduced and irreducible affine scheme over an
algebraically closed field $K$ of arbitrary characteristic. Then,
there exists an open subscheme $Y_1\subset Y$ (namely, the set of
unramified points for $\pi$) such that
$res(\pi):{\pi^{-1}(Y_1)}\rightarrow Y_1$ is an \'etale morphism
(see \cite{mum}, Ch.III.10; here, $res(\pi)$ denotes the
restriction of $\pi$.).

Let $Y,Z$ be affine, say $Y=Spec\,A, Z=Spec\,B$, where $A,B$ are
finitely generated $K$-algebras. Further, let $A,B$ be integral
domains, and $B$ an integral extension of $A$ which is finitely
generated as an $A$ module (so that $\pi:Z\longrightarrow Y$ is a
finite morphism).





Suppose that $\kappa(B)$ is a finite separable extension of
$\kappa(A)$. Then, there exists an $s\in B$ such that
$\kappa(B)=\kappa(A)[s]$. Indeed if ${{b}\over{b'}}$ is a
primitive element for the extension $\kappa(B)$ of $\kappa(A)$
with $b,b'\in B$, then there exists a $\lambda\in\kappa(A)$ such
that  $\kappa(B)=\kappa(A)[b+\lambda b']$. To see this, observe
that since there are only finitely many intermediate fields
between $\kappa(A)$ and $\kappa(B)$ (while $\kappa(A)$ is
infinite), we have that not all $\kappa(A)[b+\lambda
b'],\lambda\in \kappa(A)$ can be distinct. Hence, for some
$\lambda,\mu\in \kappa(A),\lambda\not=\mu,b+\mu b'$ is in
$\kappa(A)[b+\lambda b']$. From this it follows that
${{b}\over{b'}}$ is in $\kappa(A)[b+\lambda b']$.

\begin{lem}\label{above'}  With the above notation, let $Y=Spec\,A$ 
and let $Z=Spec\,B$.
 Consider the finite morphism
$\pi:Z\longrightarrow Y$ induced by the inclusion $A\subset B$. As
above, let $b$ be a primitive element such that
$\kappa(B)=\kappa(A)[b]$. Let $a\in A$ be such that
$A[1/a][b]=B[1/a]$. Further, let the discriminant of $b$ be
invertible in $A[1/a]$. Then $res(\pi): Spec\,B[1/a]\rightarrow
Spec\,A[1/a]$ is \'etale.

\ni (Here, $res(\pi)$ denotes the restriction of $\pi$.)

\end{lem}

\begin{proof}
Let $U=Spec\,(A[1/a]),V=\pi^{-1}(U)=Spec\,(B[1/a])$. Since
$B[1/a]=A[1/a][b]$ (by Hypothesis), we obtain that $B[1/a]$ is a
free $A[1/a]$-module with basis

\ni $\{1,b,\cdots, b^{N-1}\}$ (here, $N=[\kappa(B):\kappa(A)]$).
Further, by Hypothesis, discriminant$\,(b)$ is invertible in
$A[1/a]$. Hence we obtain that $res(\pi):
Spec\,(B[1/a])\rightarrow Spec\,(A[1/a])$ is \'etale.
\end{proof}

\begin{prop}\label{regular}

The variety $Spec\,S$ is regular in codimension $1$.

\end{prop}

\begin{proof} Taking $A=R^{O_n(K)},B=S$, as seen in the proof of
Proposition \ref{birational}, we have that $\kappa(B)$ is a
quadratic extension of $\kappa(A)$. Further, the generators $u_I$
(of the $A$-algebra $B=A[u(I),I\in I(n,m)]$) satisfy the relation
$u(I)^2-p(I,I)=0$ over $A$ (note that $A=K[p(I,J), I,J\in
I(r,m),r\le n]$ (cf. Theorem \ref{fund'})). Hence, $\kappa(B)$ is
also separable over $\kappa(A)$, since $char(K)\neq 2$. Now we
take $Y=Spec\,A$, $Z=Spec\,B$, and $\pi:Z\rightarrow Y$ the
morphism induced by the inclusion $R^{O_n(K)}\subseteq S$, in
Lemma \ref{above'}. Fixing a particular $I\in I(n,m)$, we have
that the relation $u(I)u(J)=p(I,J)$ implies that
$u(J)=p(I,J)u(I)/p(I,I)$ and hence $u(J)\in S[1/p(I,I)]$ for all
$J$. In particular $S[1/p(I,I)]$ is a free
$R^{O_n(K)}[1/p(I,I)]$-module with basis $\{1,u(I)\}$. Also the
discriminant $\delta(u_I)$ is seen to be $4p(I,I)$ which is
invertible in $R^{O_n(K)}[1/p(I,I)]$. Hence, taking
$a=p(I,I),b=u(I)$, the hypothesis of Lemma \ref{above'} are
satisfied. Hence $$res(\pi):{\underset{I\in
I(n,m)}\cup}\,Spec\,S[1/p(I,I)]\rightarrow {\underset{I\in
I(n,m)}\cup}\,Spec\,R^{O_n(K)}[1/p(I,I)]$$ is \'etale. Let
$\mathfrak{a}$ be the ideal generated by $\{p(I,I), I\in I(n,m)\}$.
Note that in view of the relations $$u(I)u(J)=p(I,J), I,J \in
I(n,m)$$ (cf. Lemma \ref{formula}), we have
$$p(I,J)^2=p(I,I)p(J,J)$$ Thus we have, $p(I,J)\in \sqrt{\mathfrak{a}},
\forall I,J \in I(n,m)$. Let $Y_0=V(\sqrt{\mathfrak{a}})$. Then $Y_0$
is simply $D_{n}(Sym\,M_m)$. Also, $Y$ being $D_{n+1}(Sym\,M_m)$
(cf. Theorem \ref{fund}), we have $Y_0$ is the singular locus of
$Y$ (cf. Theorem \ref{sing}). Hence, $codim_YY_0\ge 2$ (since $Y$
is normal).  Now the branch locus $Y_b$ of $\pi:Z\longrightarrow
Y$ is contained in $Y_0$; hence
$$codim_YY_b\ge 2\leqno{(*)}$$ Denote
$$Y_e:=\{\texttt{unramified points for }\pi\}$$ We have, $res(\pi):\pi^{-1}(Y_e)\rightarrow
Y_e$ is \'etale, and $Y_b=Y\,\setminus\,Y_e$. Denote
$$Z_e:=\pi^{-1}(Y_e),\ Z_b:=\pi^{-1}(Y_b)$$ Denoting by
$Z_{\texttt{sing}}$ the singular locus of $Z$, we have (cf. (*))
$$codim_Z(Z_{\texttt{sing}}\cap Z_b)\,(\ge codim_ZZ_b)\ge 2\leqno{(**)}$$ On the
other hand, $res(\pi):Z_e\rightarrow Y_e$ being \'etale, we have,
$$\pi^{-1}(Y_{\texttt{sing}}\cap
Y_e)=Z_{\texttt{sing}}\cap Z_e$$ Hence we obtain,
$$codim_Z(Z_{\texttt{sing}}\cap Z_e)=codim_Y(Y_{\texttt{sing}}\cap
Y_e)\ge codim_Y(Y_{\texttt{sing}})\ge 2\leqno{(***)}$$ (**) and
(***) imply that $codim_Z(Z_{\texttt{sing}})\ge 2$, and the result
follows.
\end{proof}

The above Proposition together with Theorem \ref{main} implies the
following
\begin{prop}\label{normal'}
$Spec\,S$ is normal.
\end{prop}
The result follows from Serre's criterion for normality: $Spec\,A$
is normal if and only if $A$ has $S_2$ and $R_1$.

\begin{thm}\label{equal} The inclusion $S\subseteq R^{SO_n(K)}$ is an equality, i.e., $R^{SO_n(K)}=S$.
\end{thm}
\begin{proof}
Propositions \ref{birational}, \ref{normal'} imply (in view of
Zariski Main Theorem (cf. \cite{mum}, ch III.9)) that the morphism
$q:Spec\,R^{SO_n(K)}\longrightarrow Spec\,S$ is in fact an
isomorphism. The result follows from this.
\end{proof}

Combining the above Theorem with Theorems \ref{present'},
\ref{main} we obtain the following theorems.

\begin{thm}\label{present''}\textbf{A presentation for $R^G$:}
\begin{enumerate}\item (\textbf{First fundamental Theorem}) The
$K$-algebra $R^G(=S)$ is generated by $\{p(A,B),u(I), (A,B)\in
H_p,I\in H_u\}$. \item (\textbf{Second fundamental Theorem}) The
ideal of relations among the generators $\{p(A,B),u(I)\}$ is
generated by  relations (1)-(3) of Proposition \ref{relations}.
\end{enumerate}
\end{thm}
\begin{thm}\label{main'} $R^{SO_n(K)}$ is Cohen-Macaulay.
\end{thm}
\section{Application to moduli problem}\label{moduli}
\setcounter{subsubsection}{0} In this section, using Theorem
\ref{main'}, we give a characteristic-free
proof of the Cohen-Macaulayness of the moduli space
$\modulitwo$ of equivalence classes of semi-stable rank $2$,
degree $0$ vector bundles on a smooth projective curve of genus $>
2$ by relating it to $K[X]^{SO_3(K)}$. It is known \cite[\S7,~Theorem~3]{nr} 
that $\mathcal{M}_2$ is smooth when the genus is~$2$.  

Assume for the moment that the characteristic of the field $K$ is zero.
Consider the moduli space $\mathcal{M}_n$ of
equivalence classes of semi-stable, rank $n$, degree $0$ vector
bundles on a smooth projective curve $C$ of genus $m>2$. Let $V$
be the trivial vector bundle on $C$. The automorphism group of $V$
($\cong H^0 (C, {\rm Aut}~V)$) can be identified with $GL_n(K)$.
The tangent space at $V$ of the versal deformation space (cf.
\cite{sch}) of $V$ is $H^1 (C, {\rm End}~V$) and hence it can be
identified with $m$ copies of the space $M_n(K)$ of $n \times n$
matrices, $V$ being identified with the
``origin".
The canonical action of ${\rm Aut}~V$ on this tangent
space gets identified with the diagonal adjoint action of
$GL_n(K)$ on $m$ copies of $M_n(K)$. Now the moduli space $\moduli_n$
is a GIT
quotient $Z/\!\!/ H$, for a suitable $Z$ and $H$ a projective linear group
of suitable rank. The versal deformation space of $V$ gets embedded in
$Z$ and by ``Luna slice" type of arguments (cf. \cite{git},
Appendix to chapter 1, D), the analytic local ring of $\moduli_n$
at $V$ gets identified with the analytic local ring of $H^1 (C,
{\rm End}~V)/\!\!/PGL_n(K)$ at the
``origin".

Suppose now that $n=2$.  If $V$ is a stable vector bundle,  then the
point in  $\modulitwo$ that it defines is smooth.   If $V$ is not stable,
then the point in $\modulitwo$ that it defines can be
represented by $L_1\oplus L_2$, where $L_1$ and $L_2$ are line bundles of
degree zero. If $L_1\cong L_2$, then $End\,(V)\cong M_2(K)$, and
the considerations are the same as for the case when $V$ is
trivial. If $L_1$ is not isomorphic to $L_2$, then $Aut\,V$ is the
torus $T:=\mathbb{G}_m\times \mathbb{G}_m$, and the analytic
local ring of $\modulitwo$ at $V$ is isomorphic to the analytic
local ring of $H^1(C,End\,V)/\!\!/T$ at the origin (in fact
the action of $T$ is trivial). This is certainly Cohen-Macaulay.

Although we are only interested in the rank $2$ case,  let us remark that
considerations of the previous paragraph
hold when the rank of $V$ is arbitrary, and
the analytical local ring of $\moduli_n$ is isomorphic to
$Z/\!\!/G$, where $Z=\prod\,Z_i, G=\prod\,G_i$, where $Z_i$ is $g$
copies of $M_{n_{i}}(K)$ and $G_i=PGL_{n_{i}}(K)$.

Although we have assumed that the characteristic of the field to be
zero,   the above considerations remain valid in positive characteristic
but with certain restrictions;
for example, if $n=2$, then the characteristic should not be $2$ or $3$.

Thus, in order to prove that $\modulitwo$ is Cohen-Macaulay, it suffices
to show that the point corresponding to the trivial bundle is so.
Further, since the analytic local ring at the point corresponding to the
trivial bundle gets identified with the completion of the local ring
at the origin of $M_2(K)^{\oplus m}/\!\!/PGL_2(K)$,   it suffices
to prove that $M_2(K)^{\oplus m}/\!\!/PGL_2(K)$ is Cohen-Macaulay.

\medskip Let now $n=2$, $M_2:= M_{2}(K),M_2^0:=sl_2(K) $. Let
$$Z=M_2\oplus \cdots \oplus M_2 (g\ {\rm{copies}}),\,
Z_0=M_2^0\oplus \cdots \oplus M_2^0 (g\ {\rm{copies}})$$

\medskip Let $A=K[Z],A_0=K[Z_0]$. Let $G=SL_2(K)$. Consider the diagonal
action of $G$ on $Z,Z_0$ induced by the action of $G$ on $M_2,
M_2^0$ by conjugation respectively.

From the above discussion, we have that the completion of the
local ring at the point in $\mathcal{M}_2$ corresponding to the
trivial rank $2$ vector bundle is isomorphic to the completion of
$A^G$ at the point which is the image of the origin (in
$Z(={\mathbb{A}}^{4g})$) under $Z\rightarrow \,{\rm{Spec}}A^G$. On the
other hand, we have, $A^G=A_0^G[x_1,\cdots ,x_g], x_i$'s being
indeterminates (since, $M_2\cong M^0_2\oplus K$); further, we have
that the adjoint action of $SL_2(K)$ on $sl_2(k)$ is isomorphic to
the natural representation of $SO_3(K)$ on $K^3$ (note that the
Lie algebras $sl_2(K)$ and $so_3(K)$ are isomorphic). Hence the
ring $A_0^G$ gets identified with $K[X]^{SO_3(K)}$, $X$ being
$V\oplus \cdots \oplus V (g\ {\rm{copies}}),\,V=K^3$. Hence we
obtain
\begin{thm}\label{mod}
\begin{enumerate}
\item The GIT quotients $Z\fs G,Z_0\fs G$ are Cohen-Macaulay.
\item The moduli space $\mathcal{M}_2$ is Cohen-Macaulay. \item We
have standard monomial bases for the co-ordinate rings of $Z\fs G$, $Z_0\fs G$.
\item We have first \& second fundamental theorems for
the co-ordinate rings of $Z\fs G$ and $Z_0\fs G$,
i.e., algebra generators, and generators for the
ideal of relations among the generators.\hfill$\Box$
\end{enumerate}
\end{thm}

\begin{rem}
The results in Theorem \ref{mod} being characteristic-free, we may
deduce from Theorem \ref{mod} that the moduli space
$\mathcal{M}_2$ behaves well under specializations; for instance,
if the curve $C$ is defined over $\mathbb{Z}$, and if
$\mathcal{M}_2({\mathbb{Z}})$ is the corresponding moduli space,
then for any algebraically closed field $K$ of characteristic
$\not= 2,3$, the base change of $\mathcal{M}_2({\mathbb{Z}})$ by
$K$ gives the moduli space $\mathcal{M}_2(K)$ over $K$.
\end{rem}
\begin{rem}
This section is motivated by \cite{m-r1}. In loc.cit, the
Cohen-Macaulayness for $Z\fs G,Z_0\fs G,\mathcal{M}_2$ are deduced
by proving the Frobenius-split properties for these spaces.
\end{rem}

\section{Results for the adjoint action of $SL_2(K)$}\label{sl2}
\setcounter{subsubsection}{0}  Consider $G=SL_2(K), ch\,K\not=
2$. Let
$$Z=\underset{m\text{ copies}}{\underbrace{sl_2(K)\oplus \cdots
\oplus sl_2(K)}}=Spec\,R,\texttt{ say}$$  Using the results of \S
\ref{cohen}, \S \ref{normal}, we shall describe a ``standard
monomial basis" for $R^G$; the elements of this basis will be
certain monomials in $tr\,(A_iA_j)$'s, and $tr\,(A_iA_jA_k)$'s.

\vs.2cm\ni\textbf{Identification of $sl_2(K)$ and $so_(K)$:} Let
$$X=\begin{pmatrix} 0&1\\ 0&0
 \end{pmatrix},H=\begin{pmatrix} 1&0\\ 0&-1
 \end{pmatrix},Y=\begin{pmatrix} 0&0\\ 1&0
 \end{pmatrix}$$ be the Chevalley basis of $sl_2(K)$.

 \vs.2cm\ni Let $\langle,\rangle$ be a symmetric non-degenerate bilinear form
  on $V=K^3$. Taking the matrix of the form $\langle,\rangle$ to be
 $J=\texttt{ anti-diagonal}(1,1,1)$, we have,
 $$\begin{gathered}SO_3(K)=\{A\in SL_3(K)\,|\,J^{-1}(^{t}A)^{-1}J=A\}\\
 so_3(K)=
 \{A\in sl_3(K)\,|\,J^{-1}(^{t}A)J=-A\}
 \end{gathered}$$
The Chevalley basis of $so_3(K)$ is given by $$X'=\begin{pmatrix}
0&{\sqrt{2}}&0\\ 0&0&-{\sqrt{2}}\\0&0&0
 \end{pmatrix},H'=\begin{pmatrix} 2&0&0\\ 0&0&0\\ 0&0&-2
 \end{pmatrix},Y'=\begin{pmatrix} 0&0&0\\ {\sqrt{2}}&0&0\\
 0&-{\sqrt{2}}&0
 \end{pmatrix}$$
The map $X\mapsto X',H\mapsto H',Y\mapsto Y'$ gives an isomorphism
 $sl_2(K)\cong so_3(K)$ of Lie algebras. Further, the map $$\theta: sl_2(K)\rightarrow K^3,
 \begin{pmatrix}
a&b\\c&-a
  \end{pmatrix}\mapsto
  (\frac{b}{\sqrt{2}},-a,\frac{c}{\sqrt{2}})$$ identifies the
  adjoint action of $SL_2(K)$ on $sl_2(K)$ with the natural action
  of $SO_3(K)$ on $K^3$; and the induced map $$\theta_m:
  Z\longrightarrow \underset{m\text{
copies}}{\underbrace{V\oplus \cdots \oplus V}}\,  ,\, V=K^3$$
   identifies the diagonal (adjoint) action of $SL_2(K)$ on $Z$ with the diagonal (adjoint) action of
  $SO_3(K)$ on ${\underbrace{V\oplus \cdots \oplus V}}$.

 \begin{lem}\label{trace} Let $z\in Z $, say $z=(A_1,\cdots,A_m)$.
Let $\theta_m(z)= (u_1,\cdots,u_m)$. We have
\begin{enumerate}
\item tr$\,(A_iA_j)=2\langle u_i,u_j\rangle$. \item Let $I\in
I(3,m)$, say, $I=(i,j,k)$. Then $tr\,(A_iA_jA_k)=2u(I)$.
\end{enumerate}

\ni (here, $u(I)$ is as in \S \ref{alg})
  \end{lem}
The proof is an easy verification.

\vs.2cm\ni\textbf{Notation:} In the sequel we shall denote
$$U(i,j):=tr\,(A_iA_j),i\ge j,\ U(i,j,k):=tr\,(A_iA_jA_k),(i,j,k)\in I(3,m)$$

\begin{rem}
Note that if $i,j,k$ are not distinct, then $tr\,(A_iA_jA_k)=0$;
for, say, $i=j$, then $A_i$ being a $2\times 2$ traceless matrix,
we have (in view of Cayley-Hamilton theorem), $A_i^2=-|A_i|$, and
hence
$tr\,(A_i^2A_k)=-|A_i|tr(\,A_k)=0$ (since $A_k\in sl_2(K)$).
\end{rem}

 In view of the identification given by $\theta_m$, we obtain
 (cf. Theorems \ref{equal}, \ref{capMain})
  a ``standard monomial basis" for $S:=R^{SL_2(K)}$. By Theorem \ref{capMain}, monomials
  $$p({\underline{\alpha}})p({\underline{A}},{\underline{B}})u({\underline{I}})
  $$ where
  $$\begin{gathered}
p({\underline{\alpha}}):=p(\alpha_1,\alpha_2)\cdots
p(\alpha_{2r-1},\alpha_{2r}), \mathrm{\ for\ some\ }r,\,
\alpha_i\in [1,m]\\
p({\underline{A}},{\underline{B}}):=p(A_1,B_1)\cdots
  p(A_s,B_s),  \mathrm{\ for\ some\ }s,\, A_i,B_i\in I(2,m)\\
u({\underline{I}}):=u(I_1)\cdots u(I_t),  \mathrm{\ for\ some\
}t,\, I_\ell\in I(3,m)\\
\alpha_1\ge\cdots\ge\alpha_{2r}\ge A_1\ge B_1\ge A_2\ge\cdots\ge
B_s\ge I_1\ge\cdots\ge I_\ell
  \end{gathered}$$
give a basis for $R^{SO_3(K)}$ (here, $[1,m]$ denotes the set
$\{1,\cdots,m\}$).

\ni We shall refer to
$p({\underline{\alpha}})p({\underline{A}},{\underline{B}})u({\underline{I}})
  $ as a {\em standard monomial of multi-degree} $(r,s,t)$.

\ni Denote $$\mathcal{M}_{r,s,t}:=\{\texttt{standard monomials of
multi-degree }(r,s,t)\}$$

\vs.2cm\ni\textbf{Standard monomials of Type I,II,III:} Let
notation be as above. We shall refer to
$$\begin{gathered}
p(\alpha_1,\alpha_2)\cdots p(\alpha_{2r-1},\alpha_{2r}), \,
\alpha_1\ge\cdots\ge\alpha_{2r}\\
p(A_1,B_1)\cdots
  p(A_s,B_s),  \,A_1\ge B_1\ge A_2\ge\cdots\ge
B_s\\
u(I_1)\cdots u(I_t),\,I_1\ge\cdots\ge I_\ell
  \end{gathered}$$ as {\em standard monomials of Type I,II,III} (and of degree $r,s,t$)
  respectively.

  We now define three types of standard monomials in $U(i,j)(=tr\,(A_iA_j)),i\ge j,\,U(i,j,k)
  (=tr\,(A_iA_jA_k)), (i,j,k)\in I(3,m)$ analogous to the above
  three types. Type I and Type III are direct extensions to
  traces:
  \begin{defn}
A monomial of the form $$U(\alpha_1,\alpha_2)\cdots
U(\alpha_{2r-1},\alpha_{2r}), \, \alpha_1\ge\cdots\ge\alpha_{2r}$$
will be called a Type I standard monomial.
  \end{defn}
  \begin{defn} A monomial of the form $$U(I_1)\cdots U(I_t),\,I_1\ge\cdots\ge I_\ell$$
  will be called a Type III standard monomial.
 \end{defn}
  For defining Type II standard monomials, we first define a
  bijection between

  \ni $\{p(A,B),A,B\in I(2,m),A\ge B\}$ and
  $\{U(j,i)U(l,k),j\ge i,l\ge k,i\not\ge l\}$.
Note that given $U(j,i)U(l,k),j\ge i,l\ge k$,
  we may suppose that $j$ is the greatest among $\{i,j,k,l\}$
  (if $l$ is the greatest, then we may write $U(j,i)U(l,k)$ as
  $U(l,k)U(j,i)$); then the latter set is simply
  $$\{\texttt{all non-standard Type I, degree $2$ monomials in the }
  U(a,b)\texttt{'s}\}$$   Let
  us denote the two sets by $\mathcal{A},\mathcal{B}$
  respectively. Let $p(A,B)\in \mathcal{A}$, say, $A=(a,b), b>a,B=(c,d),d>c$. Define $\omega:\mathcal{A}\rightarrow\mathcal{B}$
  as follows: $$\omega(p(A,B))=\begin{cases} U(b,c)U(d,a), \mathrm{\ if\ }d\ge
  a\\ U(b,d)U(a,c),\mathrm{\ if\ } d<a
 \end{cases}$$ Given a standard monomial $p(A_1,B_1)\cdots
  p(A_s,B_s),A_i,B_i\in I(2,m)$, we shall define $$\omega(p(A_1,B_1)\cdots
  p(A_s,B_s)):=\omega(p(A_1,B_1))\cdots
 \omega(p(A_s,B_s))$$
 \begin{defn}
A monomial $$U(j_1,i_1)U(l_1,k_1)\cdots
U(j_s,i_s)U(l_s,k_s),\,j_t>i_t,l_t>k_t,\,1\le t\le s$$ is called a
Type II standard monomial if it equals $\omega(p(A_1,B_1)\cdots
  p(A_s,B_s))$ for some (Type II) standard monomial $p(A_1,B_1)\cdots
  p(A_s,B_s)$.
 \end{defn}
 \begin{rem}
Note in particular that for $s=1$, the Type II standard monomials
(in the $U(j,i)$'s) are precisely the non-standard Type I, degree
$2$ monomials in the $U(a,b)$'s.
 \end{rem}
 Let us extend $\omega$ to  $\mathcal{M}_{r,s,t}$, the definition of $ \omega(\mathcal{F})$, for $\mathcal{F}$ a Type I or
 III standard monomial being obvious, namely, $$\begin{gathered} \omega(p(\alpha_1,\alpha_2)\cdots
 p(\alpha_{2r-1},\alpha_{2r}))=U(\alpha_1,\alpha_2)\cdots
 U(\alpha_{2r-1},\alpha_{2r})\\ u(I_1)\cdots u(I_t)=U(I_1)\cdots U(I_t)
 \end{gathered}$$
 Define $$\omega(p({\underline{\alpha}})p({\underline{A}},{\underline{B}})u({\underline{I}}))=
 \omega(p({\underline{\alpha}}))
 \omega(p({\underline{A}},{\underline{B}}))\omega(u({\underline{I}})),\
 p({\underline{\alpha}})p({\underline{A}},{\underline{B}})u({\underline{I}})\in \mathcal{M}_{r,s,t}$$
 \begin{defn}
A monomial in the $U(j,i)$'s, $j\ge i$ and $U(I)$'s, $I\in I(3,m)$
of the form
$$\omega(p({\underline{\alpha}})p({\underline{A}},{\underline{B}})u({\underline{I}})),\
p({\underline{\alpha}})p({\underline{A}},{\underline{B}})u({\underline{I}})\in
\mathcal{M}_{r,s,t}$$ will be called a standard monomial of type
$(r,s,t)$.
 \end{defn}
 \vs.2cm\ni \textbf{Notation:} We shall denote $$
\mathcal{N}_{r,s,t}=\omega(\mathcal{M}_{r,s,t})
 $$
\begin{thm}\label{standard} Let $Z=\underset{m\text{ copies}}{\underbrace{sl_2(K)\oplus \cdots
\oplus sl_2(K)}}=Spec\,R$, say, where, $m>3$. Standard monomials
(in the traces) of type $(r,s,t)$,  $r,s,t$ being non-negative
integers, form a basis for $R^{SL_2(K)}$ for the adjoint action of
$SL_2(K)$ on $Z$
\end{thm}
\begin{proof}
Denote $S:=R^{SL_2(K)}$. Write
$$S={\underset{(r,s,t)}{\oplus}}{S_{r,s,t}}$$ where $r,s,t$ are positive
integers, and $S_{r,s,t}$ is the $K$-span of
$\mathcal{M}_{r,s,t}$. In fact, we have (in view of linear
independence of standard monomials (cf. Theorem \ref{capMain})) that
$\mathcal{M}_{r,s,t}$ is a basis for $S_{r,s,t}$. Using the
bijection $\omega$, we shall show that $\mathcal{N}_{r,s,t}$ is
also a basis for $S_{r,s,t}$. Clearly this requires a proof only
in the case $s\not= 0$. Let $N_{r,s,t}=\#\mathcal{M}_{r,s,t}$. The
relations (cf. Lemma \ref{trace}; note that $p(i,j)=\langle
u_i,u_j\rangle$)
$$\begin{gathered}
 p(i,j)=\frac{1}{2}U(i,j)\\
 p(A,B)=\frac{1}{4}(U(a,c)U(b,d)-U(b,c)U(a,d))\\
 u(I)=\frac{1}{2}U(I)
 \end{gathered}$$ give raise to the transition matrix, say, $M$. Then it is easy to see that
  for a suitable indexing of the elements of $\mathcal{M}_{r,s,t}$, and
 $\mathcal{N}_{r,s,t}$, the matrix $M$ takes the upper triangular form
 with the diagonal entries being non-zero. To be very precise, to
 $p({\underline{\alpha}})p({\underline{A}},{\underline{B}})u({\underline{I}})
  \in \mathcal{M}_{r,s,t}$, we associate a
  $(4s+2r+3t)$-tuple
  $n({\underline{\alpha}},{\underline{A}},{\underline{B}},{\underline{I}})$ as follows.
   Let $$\begin{gathered}
p({\underline{\alpha}})=p(\alpha_1,\alpha_2)\cdots
p(\alpha_{2r-1},\alpha_{2r}), \,
\alpha_1\ge\cdots\ge\alpha_{2r}\\
p({\underline{A}},{\underline{B}})=p(A_1,B_1)\cdots
  p(A_s,B_s),  \,A_1\ge B_1\ge A_2\ge\cdots\ge
B_s\\
u({\underline{I}})=u(I_1)\cdots u(I_t),\,I_1\ge\cdots\ge I_\ell\\
A_i=(a_{i1},a_{i2}),a_{i1}<a_{i2},\,
  B_i=(b_{i1},b_{i2}),b_{i1}<b_{i2},\,A_i\ge B_i\ 1\le i\le s
  \end{gathered}$$
   Set
  $$n({\underline{\alpha}},{\underline{A}},{\underline{B}},{\underline{I}})=
  (a_{12},a_{11},b_{12},b_{11},a_{22},a_{21}\cdots,b_{s2},b_{s1},{\underline{\alpha}},
  {\underline{I}})$$ where
  ${\underline{\alpha}}=(\alpha_1,\cdots,\alpha_{2r})$, and
  ${\underline{I}}=(I_1,\cdots,I_t )$. We take an indexing on
  $\mathcal{M}_{r,s,t}$ induced by the lexicographic order on the
  $n({\underline{\alpha}},{\underline{A}},{\underline{B}},{\underline{I}})$'s,
  and take the induced indexing on $\mathcal{N}_{r,s,t}$ (via the bijection
  $\omega$). With respect to these indexings on
 $\mathcal{M}_{r,s,t},\mathcal{N}_{r,s,t}$, it is easily seen that the transition matrix
 $M$ is upper triangular
 with the diagonal entries being non-zero. It follows that
 $\mathcal{N}_{r,s,t}$ is a basis for $S_{r,s,t}$.
\end{proof}
\begin{rem}
Note that the standard monomial basis (in the traces) as given by
Theorem \ref{standard} is characteristic-free. Also, note that
Theorem \ref{standard} recovers the result of \cite{pro}, Theorem
3.4(a), for the case of $SL_2(K)$, namely,
$tr\,(A_iA_j),tr\,(A_iA_jA_k),i,j,k\in [1,m]$ generate $R^G$ as a
$K$-algebra (in a characteristic-free way).
\end{rem}

\end{document}